\pdfoutput=1
\documentclass[toc]{mathprint}
\usepackage{rutarmacros}
\usepackage{macros}
\title{Tangents and slices of self-affine carpets}
\author
  {Antti Käenmäki}
  {University of Eastern Finland,
    Department of Physics and Mathematics,
    P.O.\ Box 111,
    FI-80101 Joensuu,
    Finland}
  {antti@kaenmaki.net}

\author
  {Alex Rutar}
  {Department of Mathematics and Statistics,
      University of Jyväskylä,
      P.O.\ Box 35 (MaD),
      FI-40014 University of Jyväskylä,
  Finland}
  {alex@rutar.org}

\begin{document}

\begin{abstract}
    We study the fine scaling properties of planar self-affine carpets.
    For Gatzouras--Lalley carpets, we give a precise formula for maximal Hausdorff dimension of a tangent in terms of the Hausdorff dimension of the projection and the Assouad dimension of the corresponding vertical slice.
    Using regularity properties for the Assouad dimension of non-autonomous self-similar sets, this implies that the set of points with tangents that are as large as possible has full Hausdorff measure, at the critical exponent.
    On the other hand, we give an explicit example of a Barański carpet for which the Hausdorff dimension of the set of points for which there exists a maximal tangent has Hausdorff dimension strictly less than the Hausdorff dimension of the original carpet.
\end{abstract}

\section{Introduction}

A classical problem in geometric measure theory is the following: given a subset $E\subset\R^d$, what can be said about the structure of the set of tangents at points in $E$?
If $E$ has positive and finite $s$-dimensional Hausdorff measure, then classical density theorems (see, for instance, \cite[Theorem~6.2]{zbl:0819.28004}) imply that almost every point has a tangent with positive Hausdorff $s$-measure.
However, sets which are not Ahlfors regular can have points with tangents which are much larger than expected.
In general, the maximal possible Hausdorff dimension of a tangent is given by the \emph{Assouad dimension} of the set $E$ \cite{zbl:1154.37322,doi:10.1093/imrn/rnw336}; this value is attained by a \emph{weak tangent} (where the location of ``zooming in'' is allowed to vary depending on the scale), but not necessarily by an actual tangent \cite{zbl:1321.54059}.

Moreover, for many well-studied fractal sets, the Assouad dimension and the Hausdorff dimension can differ, so the classical information concerning the existence of large tangents cannot reach the threshold of the Assouad dimension.
Motivated by this phenomenon, in \cite{arxiv:2309.11971}, the authors study the structure of the set of tangents for sets satisfying various weak forms of dynamical invariance (informally, that the sets contain small potentially highly distorted copies at all scales and distortions).
For general attractors of bi-Lipschitz iterated function systems, it is shown in \cite{arxiv:2309.11971} that there necessarily is at least one tangent with Hausdorff dimension which attains the Assouad dimension, and (essentially) overlapping self-conformal sets, the authors prove in fact that there the set of points with tangents of maximal Hausdorff dimension is a full dimension subset.

\begin{figure}[t]
    \centering
    \begin{subcaptionblock}{.47\textwidth}
        \centering
        \includegraphics[width=0.7\textwidth]{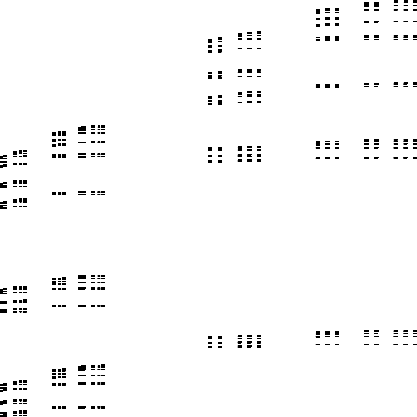}
        \caption{Gatzouras--Lalley}
    \end{subcaptionblock}%
    \begin{subcaptionblock}{.47\textwidth}
        \centering
        \includegraphics[width=0.7\textwidth]{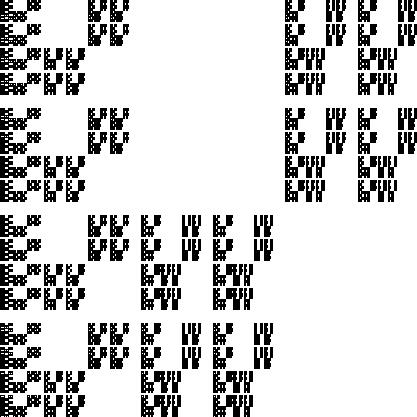}
        \caption{Barański}
    \end{subcaptionblock}%
    \caption{Some self-affine carpets, which are attractors of the iterated function systems depicted in \cref{f:carpets}.}
    \label{f:attractors}
\end{figure}

However, the gap between these two classes of sets is quite large.
For instance, one might hope that in the presence of reasonably well-behaved dynamics, then most points will have a tangent which is as large as possible.
In order to understand this problem more generally, in this paper, we study the question of the fine structure of tangents for a particular family of \emph{self-affine sets} (such as those depicted in \cref{f:attractors}).

Our results show that in fact whether or not there are many maximal tangents depends on the geometry of the particular set under consideration.
Within our class of sets, when there is only one direction of maximal contraction, the majority of points have maximal tangents which are as large as possible.
Despite this, there are still many points for which all tangents are much smaller than expected.
This is explained in \cref{it:gl-results}.
On the other hand, we also demonstrate that there exist self-affine sets for which the majority of points have no large tangents (see \cref{it:b-results}).
More generally, we obtain precise results relating the dimensions of tangents to the Assouad dimension of appropriate slices of the corresponding self-affine set.

In the next section, we introduce the relevant definitions to make the above definitions precise.
We then state our main results in \cref{ss:results}.

\subsection{Assouad dimensions}
Fix a compact set $K\subset\R^d$.
We say that a compact set $F\subset B(0,1)$ is a \defn{weak tangent} of $K\subset\R^d$ if it is a Hausdorff limit of successive magnifications of the set $K$.
We denote the set of weak tangents of $K$ by $\Tan(K)$.
More strongly, $F$ is a \defn{tangent} of $K$ at $x$ if it is a Hausdorff limit of successive magnifications of $K$ upon the point $x$.
We denote the set of tangents of $K$ at $x$ by $\Tan(K,x)$.
We refer the reader to \cref{ss:tan-def} for precise definitions.

Closely related to the notion of a weak tangent is the \defn{Assouad dimension} of $K$, introduced in \cite{zbl:0396.46035}, which is the dimensional quantity
\begin{align*}
    \dimA K=\inf\Bigl\{s:\exists C>0\,&\forall 0<r\leq R<1\,\forall x\in K\\*
                                      &N_r(B(x,R)\cap K)\leq C\Bigl(\frac{R}{r}\Bigr)^s\Bigr\}.
\end{align*}
Here, for a general bounded set $F$, $N_r(F)$ is the smallest number of closed balls with radius $r$ required to cover $F$.
It always holds that $\dimH K\leq \dimP K\leq\dimuB K\leq\dimA K$, where $\dimH K$, $\dimP K$, and $\dimuB K$ denote the Hausdorff, packing, and upper box dimensions respectively.
In some sense, the Assouad dimension is the largest reasonable notion of dimension which can be defined using covers.
An observation which goes back essentially to Furstenberg, but was stated explicitly in \cite{doi:10.1093/imrn/rnw336}, is that the Assouad dimension is characterized by weak tangents:
\begin{equation}\label{e:large-wt}
    \dimA K=\max\{\dimH F:F\in\Tan(K)\}.
\end{equation}

Continuing the analogy with tangents, a localized version of the Assouad dimension was recently introduced in \cite{arxiv:2309.11971}.
Given $x\in K$, the \defn{pointwise Assouad dimension} is
\begin{align*}
    \dimA(K,x)=\inf\Bigl\{s:\exists C>0\,&\exists\rho>0\,\forall 0<r\leq R<\rho\\*
                                         &N_r(B(x,R)\cap K)\leq C\Bigl(\frac{R}{r}\Bigr)^s\Bigr\}.
\end{align*}
The choice of $\rho>0$ in the definition of $\dimA(K,x)$ ensures a sensible form of bi-Lipschitz invariance: if $f\colon K\to K'$ is bi-Lipschitz, then $\dimA(K,x)=\dimA(f(K),f(x))$.
It is immediate from \cite[Proposition 2.2]{arxiv:2309.11971} and the definition of the pointwise Assouad dimension that
\begin{equation}\label{e:large-tan}
    \sup\{\dimuB F:F\in\Tan(K,x)\} \leq \dimA(K,x) \leq \dimA K.
\end{equation}
Unfortunately, the first inequality can be strict by \cite[Example~2.9]{arxiv:2309.11971}, and the second inequality can be strict for all $x\in K$ simultaneously by \cite[Example~2.8]{arxiv:2309.11971}.
On the other hand, if $K$ is very regular (for instance, Ahlfors--David regular), then $\dimA(K,x)=\dimA K$ for all $x\in K$.

Of course, there are many sets which are not Ahlfors--David regular, but which still exhibit enough regularity that one might hope for more to be true.
A natural question, motivated by the general relationship in \cref{e:large-wt}, is to understand when equality holds in \cref{e:large-tan}.
Some preliminary answers are given in \cite[Theorem~C]{arxiv:2309.11971}.
In particular, if $K$ is a self-affine set (or more generally an attractor of a bi-Lipschitz IFS), then
\begin{equation*}
    \dimA K=\max\{\dimH F:F\in\Tan(K,x)\text{ and }x\in K\}
\end{equation*}
and
\begin{equation} \label{e:pack-full}
  \dimP\{x\in K:\dimA(K,x)=\dimA K\}=\dimP K.
\end{equation}
Moreover, if $K$ is a self-similar or self-conformal set, then \cite[Theorem 2.12]{arxiv:2309.11971} shows that in fact equality holds on a very large subset:
\begin{equation} \label{e:haus-full}
    \dimH\{x\in K:\exists F\in\Tan(K,x)\text{ such that }\dimH F = \dimA K\}=\dimH K.
\end{equation}
However, there is a relatively large gap between self-conformal sets and general attractors of bi-Lipschitz IFSs.
Motivated by these results, our goal in this document is to understand to what extent the results for self-conformal sets extend to more general IFS attractors.

In the following section, we discuss our main results concerning self-affine sets and provide some answers which indicate that answers to these questions are, in general, quite subtle.

\subsection{Main results}\label{ss:results}
With these questions in mind, we now turn our attention to two specific families of affine iterated function systems in the plane: specifically, the planar self-affine carpets of \emph{Gatzouras--Lalley} \cite{zbl:0757.28011} and \emph{Barański} \cite{zbl:1116.28008}.
Note that these sets are self-affine but (except for some degenerate cases) not self-similar.
We defer precise definitions and notation to \cref{ss:carpet-def}; see \cref{f:carpets} for examples of the generating maps in these classes.
\begin{figure}[t]
    \centering
    \begin{subcaptionblock}{.47\textwidth}
        \centering
        \begin{tikzpicture}[scale=5]
    \begin{scope}[font=\tiny]
        \node[below] at (0,0) {$0$};
        \node[below] at (1,0) {$1$};
        \node[left] at (0,0) {$0$};
        \node[left] at (0,1) {$1$};
    \end{scope}
    \draw (0,0) rectangle (1,1);
    \begin{scope}[thick]
        \draw[fill=gray!10] (0.0,0.0) rectangle (1/4,1/8);
        \draw[fill=gray!10] (0.0,1/4) rectangle (1/4,1/4+1/11);
        \draw[fill=gray!10] (0.0,1/2) rectangle (1/4,1/2+1/5);

        \draw[fill=gray!10] (1/2,1/8) rectangle (1,1/6);
        \draw[fill=gray!10] (1/2,3/4) rectangle (1,1);
        \draw[fill=gray!10] (1/2,1/2+1/9) rectangle (1,1/2+1/6);
    \end{scope}
\end{tikzpicture}
        \caption{Gatzouras--Lalley}
    \end{subcaptionblock}%
    \begin{subcaptionblock}{.47\textwidth}
        \centering
        \begin{tikzpicture}[scale=5]
\begin{scope}[font=\tiny]
\node[below] at (0,0) {$0$};
\node[below] at (1,0) {$1$};
\node[left] at (0,0) {$0$};
\node[left] at (0,1) {$1$};
\end{scope}
\draw (0,0) rectangle (1,1);
\begin{scope}[thick]
\draw[fill=gray!10] (0.0,0.0) rectangle (0.2833333333333333,0.2);
\draw[fill=gray!10] (0.0,0.26666666666666666) rectangle (0.2833333333333333,0.4666666666666667);
\draw[fill=gray!10] (0.0,0.5333333333333333) rectangle (0.2833333333333333,0.7333333333333333);
\draw[fill=gray!10] (0.0,0.8) rectangle (0.2833333333333333,1.0);
\draw[fill=gray!10] (0.3458333333333333,0.0) rectangle (0.4625,0.2);
\draw[fill=gray!10] (0.3458333333333333,0.26666666666666666) rectangle (0.4625,0.4666666666666667);
\draw[fill=gray!10] (0.525,0.0) rectangle (0.6416666666666667,0.2);
\draw[fill=gray!10] (0.525,0.26666666666666666) rectangle (0.6416666666666667,0.4666666666666667);
\draw[fill=gray!10] (0.7041666666666667,0.5333333333333333) rectangle (0.8208333333333333,0.7333333333333333);
\draw[fill=gray!10] (0.7041666666666667,0.8) rectangle (0.8208333333333333,1.0);
\draw[fill=gray!10] (0.8833333333333333,0.5333333333333333) rectangle (1.0,0.7333333333333333);
\draw[fill=gray!10] (0.8833333333333333,0.8) rectangle (1.0,1.0);
\end{scope}
\begin{scope}[font=\tiny]
\node[below] at (0.14166666666666666, 0) {$\alpha_1$};
\node[below] at (0.4041666666666667, 0) {$\alpha_2$};
\node[below] at (0.5833333333333334, 0) {$\alpha_2$};
\node[below] at (0.7625, 0) {$\alpha_2$};
\node[below] at (0.9416666666666667, 0) {$\alpha_2$};
\node[left] at (0, 0.1) {$\beta$};
\node[left] at (0, 0.36666666666666664) {$\beta$};
\node[left] at (0, 0.6333333333333333) {$\beta$};
\node[left] at (0, 0.9) {$\beta$};
\end{scope}
\end{tikzpicture}
        \caption{Barański}\label{sf:bar}
    \end{subcaptionblock}%
    \caption{
        Generating maps associated with a Gatzouras--Lalley and Barański system.
        The parameters for the Barański carpet correspond to the example in \cref{c:exceptional-bar} with $\delta=1/20$.
    }
    \label{f:carpets}
\end{figure}
In the following statement, let $\eta\colon\R^2\to\R$ be the orthogonal projection onto the first coordinate axis and for $x\in\R^2$ let $\ell_x$ be the vertical line containing $x$.
Recall that $\mathcal{H}^{\dimH K}(K)>0$ for every Gatzouras--Lalley carpet $K$ and in fact $\mathcal{H}^{\dimH K}(K) = \infty$ in the non-trivial (i.e.\ non-Ahlfors-regular) case; see \cite{zbl:0757.28011}.
\begin{itheorem}\label{it:gl-results}
    Let $K$ be a Gatzouras--Lalley carpet.
    Then
    \begin{equation*}
        \mathcal{H}^{\dimH K}\bigl(\{x\in K:\dimA(K,x)\neq\dimA K\}\bigr)=0.
    \end{equation*}
    On the other hand, for any $\dimB K\leq\alpha\leq\dimA K$,
    \begin{equation*}
        \dimH\{x\in K:\dimA(K,x)=\alpha\}=\dimH K.
    \end{equation*}
    Moreover, if $\eta(K)$ satisfies the SSC, then for any $x\in K$,
    \begin{enumerate}[nl,r]
        \item\label{im:pointwise-dim} $\displaystyle\max\{\dimH F:F\in\Tan(K,x)\}=\dimB\eta(K)+\dimA \ell_x\cap K$,
        \item\label{im:tangent} $\displaystyle\dimA(K,x)=\max\{\dimB K,\dimB\eta(K)+\dimA \ell_x\cap K\}$.
    \end{enumerate}
\end{itheorem}
Of course, if $\alpha\notin[\dimB K,\dimA K]$, then $\{x\in K:\dimA(K,x)=\alpha\}=\varnothing$.
It follows immediately from \cref{it:gl-results} that the conclusion \cref{e:haus-full} extends to the class of Gatzouras--Lalley carpets and that
\begin{equation*}
  \dimA(K,x)=\max\{\dimH F:F\in\Tan(K,x)\}
\end{equation*}
if and only if $\dimA\ell_x\cap K\geq\dimB K-\dimB\eta(K)$.
Moreover, $\mathcal{H}^{\dimH K}(K)<\infty$ if and only if $K$ is Ahlfors--David regular (see \cite{zbl:0757.28011}), in which case the results are trivial.
We thus see that the majority of points, from the perspective of Hausdorff $s$-measure, have tangents with Hausdorff dimension attaining the Assouad dimension of $K$.
However, we still have an abundance of points with pointwise Assouad dimension giving any other reasonable value.

\begin{remark}
    Another way to think about \cref{it:gl-results} is as follows.
    If $\bm{p} = (p_i)_{i\in\mathcal{I}}$ is any probability vector, the corresponding \emph{self-affine measure} on a Gatzouras--Lalley carpet $K$ is the unique Borel probability measure $\mu = \mu_{\bm{p}}$ satisfying $\mu=\sum_{i\in\mathcal{I}}p_i\mu\circ T_i^{-1}$, where the maps $T_i$ are the generating maps as in \cref{s:carpet-tangents}.
    Let $\underline{\ell}$ be a column of maximal column dimension $t(\underline{\ell})$, so that $\dimA K=\dimB\eta(K)+t(\underline{\ell})$ by \cref{p:gl-dims}.

    Now if $p_i > 0$ for some index $i$ in the column $\underline{\ell}$, then $\bm{p}^{\N}$-a.e.\ $\gamma\in\mathcal{I}^{\N}$ contains arbitrarily long strings of consecutive indices $i$.
    Therefore $\dimA K_{\eta(\gamma)} = t(\underline{\ell})$, so it follows from \cref{l:square-structure} that
    \begin{equation*}
        \dimA(K, x) = \dimA K\qquad\text{for $\mu$-a.e.\ $x\in K$}.
    \end{equation*}
    In particular, this is the case for any fully supported self-affine measure on $K$.

    In contrast, proving that there are points with pointwise Assouad dimension attaining any other value requires more care.
    Such points are not typical for invariant measures, and instead we construct these points by combining very long scale blocks corresponding to a uniform and well-separated subset with short scale blocks to prescribe the pointwise Assouad dimension at all points in the subset.
    Unlike the case of maximal pointwise Assouad dimension, it is important that the subset is well-separated since we require a non-trivial upper bound on the pointwise Assouad dimension.
\end{remark}

The proof of \cref{it:gl-results} is obtained by combining \cref{t:gl-tangent-formula} and \cref{t:gl-large-Assouad}.
The dimensional results given in \cref{im:pointwise-dim} and \cref{im:tangent} exhibit a precise version of a well-known phenomenon: at small scales, properly self-affine sets and measures look like products of the projection with slices.

\begin{remark}\label{r:ssc-necessary}
    In order to obtain \cref{im:pointwise-dim} and \cref{im:tangent} of \cref{it:gl-results}, the strong separation condition in the projection is required or the pointwise Assouad dimension could be incorrect along sequences which are ``arbitrarily close together at small scales''.
    The formula holds for more general Gatzouras--Lalley carpets if one restricts attention to points where this does not happen (see \cref{d:regular}), and whenever we require non-trivial upper bounds on the pointwise Assouad dimension we must ensure that all points satisfy this assumption.

    It is easy to check that the strong separation condition holds if and only if every point of $K$ is regular.
    It is also not too difficult to construct examples of systems with non-regular points where \cref{im:pointwise-dim} and \cref{im:tangent} are false.
    Of course, these formulas may still be valid even at non-regular points.
\end{remark}

For Gatzouras--Lalley carpets with projection onto the first coordinate axis satisfying the strong separation condition, slices through $x$ are precisely attractors of a non-autonomous iterated function system corresponding to the sequence of columns containing the point $x$ (such a phenomenon was exploited in a more general setting in \cite{zbl:1561.28063}).
The proof of \cref{it:gl-results} therefore relies on the dimension theory of non-autonomous self-similar sets as studied in \cite{zbmath:08103993}.
We present the required results in \cref{sec:non-autonomous}.

However, it turns out that the fact that Gatzouras--Lalley carpets have an abundance of large tangents does not extend to the non-dominated setting.
\begin{itheorem} \label{it:b-results}
    There exists a Barański carpet $K$ such that
    \begin{equation*}
        \dimH\{x\in K:\dimA(K,x)=\dimA K\}<\dimH K.
    \end{equation*}
\end{itheorem}
The conclusion \cref{e:pack-full} is thus not valid for the Hausdorff dimension in general.
The proof of \cref{it:b-results} is given in \cref{c:exceptional-bar}, and it follows from a more general result---namely \cref{t:bar-const}---describing when Barański carpets satisfying certain separation conditions have a large number of large tangents.
The proof follows from formulas for the pointwise Assouad dimension at points which are coded by sequences which contract uniformly in one direction; see \cref{p:baranski-tans} for a precise formulation.

The key distinction between Barański carpets and Gatzouras--Lalley carpets is that Barański carpets are not \emph{dominated}; see for instance \cite[\S 2.4]{arxiv:2107.00983} for a precise definition in the general planar self-affine setting.
A natural question therefore is if this phenomenon is specific to the non-dominated setting.
\begin{question}
    Let $K$ be a non-empty dominated self-affine set, with or without overlaps.
    Does it necessarily hold that
    \begin{equation*}
        \dimH\{x\in K:\dimA(K,x)=\dimA K\}=\dimH K?
    \end{equation*}
    What about the corresponding question for tangents?
\end{question}

\subsection{Brief outline and notation}

In \cref{sec:pointwise}, we recall the general framework for tangents, weak tangents and pointwise Assouad dimension, and introduce the required language concerning of metric trees and non-autonomous self-similar.
Then in \cref{s:carpet-tangents}, we apply this machinery to Gatzouras--Lalley carpets.
Most of the work is to set up fine covering arguments using approximate squares and symbolic slices in order to prove \cref{it:gl-results}.
Finally, in \cref{sec:baranski} we modify the arguments to the Barański case, which requires additional covering arguments to handle local geometric structure resulting from the lack of domination, in particular proving \cref{it:b-results}.

Throughout, we work in $\R^2$ equipped with the usual Euclidean metric.
Given functions $f$ and $g$, we say that $f\lesssim g$ if there is a constant $C>0$ so that $f(x)\leq C g(x)$ for all $x$ in the domain of $f$ and $g$.
We write $f\approx g$ if $f\lesssim g$ and $g\lesssim f$.

\subsection{Further work}
It would also be natural to ask analogous questions for self-affine measures.
Given a Borel measure $\mu$ on $\R^d$, the \defn{Assouad dimension} of $\mu$ is
\begin{align*}
    \dimA\mu=\inf\Bigl\{s:\exists C>0\,\exists\rho>0\,&\forall x\in\supp\mu\,\forall 0<r\leq R<\rho\\*
                                       &\frac{\mu(B(x,R))}{\mu(B(x,r))}\leq C\Bigl(\frac{R}{r}\Bigr)^s\Bigr\}.
\end{align*}
Similarly, for $x\in\supp\mu$, the \emph{pointwise Assouad dimension} of $\mu$ at $x$ is
\begin{align*}
    \dimA(\mu,x)=\inf\Bigl\{s:\exists C> 0\,&\exists\rho>0\,\forall 0<r\leq R<\rho\\*
                                      &\frac{\mu(B(x,R))}{\mu(B(x,r))}\leq C\Bigl(\frac{R}{r}\Bigr)^s\Bigr\}.
\end{align*}
We refer the reader to \cite{zbl:1467.28001} for more background on the Assouad dimension of measures and to \cite{zbl:1527.28004} for more background on the pointwise Assouad dimension for measures.

The following question may be of interest.
\begin{question}
    Let $\bm{p}$ and $\bm{w}$ be probability vectors on $\mathcal{I}$, and let $\mu_{\bm{p}}$ and $\mu_{\bm{w}}$ denote the corresponding self-affine measures on the Gatzouras--Lalley carpet $K$.
    Suppose $\supp\mu_{\bm{w}}\subset\supp\mu_{\bm{p}}$.
    Is there a $\mu_{\bm{w}}$-a.e.\ constant value for the function $x\mapsto \dimA(\mu_{\bm{p}}, x)$, and if so what is the value?
\end{question}
In the special case that $\bm{p} = \bm{w}$ and $\bm{p}$ has $p_i > 0$ for all $i$, it is plausible that for $\mu_{\bm{p}}$-a.e.\ $x\in K$ that
\begin{equation*}
    \dimA(\mu_{\bm{p}},x)=\dimA\mu_{\bm{p}}.
\end{equation*}
In general, it seems that most important aspect of $\bm{w}$ is its support, rather than the precise value.
However, we do not pursue this matter any more in this paper.

\section{Tangents and pointwise Assouad dimension} \label{sec:pointwise}

\subsection{Tangents and weak tangents}\label{ss:tan-def}
To begin this section, we precisely define \emph{tangents} and \emph{weak tangents}, and establish the fundamental relationship between the dimensions of tangents and the pointwise Assouad dimension.
These results will be used to find a lower bound for the pointwise Assouad dimension of a Gatzouras-Lalley carpet by means of symbolic fibres.

Given a set $E\subset\R^d$ and $\delta>0$, we denote the \emph{open $\delta$-neighbourhood} of $E$ by
\begin{equation*}
    E^{(\delta)}=\{x\in\R^d:\exists y\in E\text{ such that }|x-y|<\delta\}.
\end{equation*}
Now given a non-empty subset $X\subset\R^d$, we let $\mathcal{K}(X)$ denote the set of non-empty compact subsets of $X$ equipped with the \defn{Hausdorff metric}
\begin{equation*}
    d_{\mathcal{H}}(K_1,K_2)=\max\{p_{\mathcal{H}}(K_1;K_2),p_{\mathcal{H}}(K_2;K_1)\}
\end{equation*}
where
\begin{equation*}
    p_{\mathcal{H}}(K_1;K_2)=\inf\{\delta>0:K_1\subset K_2^{(\delta)}\}.
\end{equation*}
If $X$ is compact, then $(\mathcal{K}(X),d_{\mathcal{H}})$ is a compact metric space itself.
We also write
\begin{equation*}
    \dist(E_1,E_2)=\inf\{|x-y|:x\in E_1, y\in E_2\}
\end{equation*}
for non-empty sets $E_1,E_2\subset\R^d$.

We say that a set $F\in\mathcal{K}(B(0,1))$ is a \defn{weak tangent} of $K\subset\R^d$ if there exists a sequence of similarity maps $(T_k)_{k=1}^\infty$ with $0\in T_k(K)$ and similarity ratios $\lambda_k$ diverging to infinity such that
\begin{equation*}
    F=\lim_{k\to\infty}T_k(K)\cap B(0,1)
\end{equation*}
in $\mathcal{K}(B(0,1))$.
We denote the set of weak tangents of $K$ by $\Tan(K)$.
A key feature of the Assouad dimension is that it is characterized by Hausdorff dimensions of weak tangents.
This result is originally from \cite[Proposition~5.7]{doi:10.1093/imrn/rnw336}.
The following version is the recent improvement \cite[Corollary B]{arxiv:2309.11971}.
\begin{proposition}\label{c:full-tan}
    If $K$ is a compact set, then
    \begin{equation*}
        \alpha\coloneqq\dimA K=\max_{F\in\Tan(K)}\dimH F.
    \end{equation*}
    Moreover, the maximizing weak tangent $F$ can be chosen so that $\mathcal{H}^{\alpha}_\infty(F) \ge 1$.
\end{proposition}
In a similar flavour, we say that $F$ is a \defn{tangent of $K$ at $x\in K$} if there exists a sequence of similarity ratios $(\lambda_k)_{k=1}^\infty$ diverging to infinity such that
\begin{equation*}
    F=\lim_{k\to\infty}\lambda_k(K-x)\cap B(0,1)
\end{equation*}
in $\mathcal{K}(B(0,1))$.
We denote the set of tangents of $K$ at $x$ by $\Tan(K,x)$.

Of course, $\Tan(K,x)\subset\Tan(K)$.
Unlike in the case for weak tangents, we require the similarities in the construction of the tangent to in fact be homotheties.
This choice is natural since, for example, a function $f\colon\R\to\R$ is differentiable at $x$ if and only if the set of tangents of the graph of $f$ at $(x,f(x))$ is the singleton $\{B(0,1)\cap\ell\}$ for some non-vertical line $\ell$ passing through the origin.
In practice, compactness of the group of orthogonal transformations in $\R^d$ means this restriction will not cause any technical difficulties.

Let us next recall from \cite{arxiv:2309.11971} some of the recent results on the connections between tangent sets and pointwise Assouad dimension.
The first result is \cite[Proposition 2.2]{arxiv:2309.11971}.

\begin{proposition}\label{p:pointwise-tangent}
    For any compact set $K\subset\R^d$ and $x\in K$, $\dimA(K,x)\geq\dimuB F$ for any $F\in\Tan(K,x)$.
\end{proposition}
The next result follows from \cite[Proposition 2.11]{arxiv:2309.11971}.
\begin{proposition}\label{p:tan-lower}
    Let $K\subseteq\R^d$ be a self-affine set.
    Then for all $x\in K$, we have $\dimA(K,x)\geq\dimuB K$.
\end{proposition}
We also recall from \cite[Proposition 3.1]{zbmath:08103993} an alternative characterization of the Assouad dimension by localized packings of balls which may have very different sizes.
Let $X$ be a bounded metric space.
If $x\in X$ and $R\in(0,1)$, then the family of all localized centred packings is
\begin{equation*}
    \operatorname{pack}(X,x,R)=\left\{\{B(x_i,r_i)\}_i:\begin{matrix}0<r_i\leq R,x_i\in X,B(x_i,r_i)\subset B(x,R),\\B(x_i,r_i)\cap B(x_j,r_j)=\varnothing\text{ for all }i\neq j\end{matrix}\right\}.
\end{equation*}
The collections here may be finite or countably infinite.
In a similar way to how box and packing dimensions are related (see, for instance, \cite[Section~2.6]{zbl:1390.28012}), we have the following ``disc-packing'' formulation of the Assouad dimension.
\begin{proposition}\label{p:Assouad-disc-packing}
    Let $X$ be a bounded metric space.
    Then
    \begin{align*}
        \dimA X=\inf\Bigl\{\alpha:\forall 0<R<1\,\forall x\in X\,&\forall\{B(x_i,r_i)\}_{i=1}^\infty\in\operatorname{pack}(X,x,R)\\*
                                               &\sum_{i=1}^\infty r_i^\alpha\lesssim_\alpha R^\alpha\Bigr\}.
    \end{align*}
\end{proposition}

\subsection{Metric trees}\label{ss:metric}
First, fix a reference set $\Omega$ and write $\mathcal{T}_0=\{\Omega\}$.
Let $\{\mathcal{T}_k\}_{k=1}^\infty$ be a sequence of countable partitions of $\Omega$ so that $\mathcal{T}_{k+1}$ is a refinement of the partition $\mathcal{T}_k$.
For each $Q\in\mathcal{T}_k$ with $k\in\N$, there is a unique \defn{parent} $\widehat{Q}\in\mathcal{T}_{k-1}$ with $Q\subset\widehat{Q}$.
Suppose that for any $\gamma_1\neq\gamma_2\in\Omega$ there is a $k\in\N$ such that there are $Q_1\neq Q_2\in\mathcal{T}_k$ so that $\gamma_1\in Q_1$ and $\gamma_2\in Q_2$.
We call such a family $\{\mathcal{T}_k\}_{k=0}^\infty$ a \defn{tree}, and write $\mathcal{T}=\bigcup_{k=0}^\infty\mathcal{T}_k$.

Now, suppose that there is a function $\rho\colon\mathcal{T}\to(0,\infty)$ which satisfies
\begin{enumerate}[nl]
    \item\label{im:rho-sep} $0<\rho(Q)<\rho(\widehat{Q})$, and
    \item\label{im:rho-finite} there is a sequence $(r_k)_{k=1}^\infty$ converging to zero from above such that $\rho(Q)\leq r_k$ for all $Q\in\mathcal{T}_k$.
\end{enumerate}
The function $\rho$ induces a metric $d$ on the space $\Omega$ by the rule
\begin{equation*}
    d(\gamma_1,\gamma_2)=\inf\{\rho(Q):Q\in\mathcal{T}\text{ and }\{\gamma_1,\gamma_2\}\subset Q\}.
\end{equation*}
In particular, $\diam(Q)=\rho(Q)$ with respect to the metric $d$.
We then refer to the data $(\Omega,\{\mathcal{T}_k\}_{k=0}^\infty,\rho)$ as a \defn{metric tree}.

We say that a subset $\mathcal{A}\subset\mathcal{T}$ is a \defn{section} if $Q_1\cap Q_2=\varnothing$ whenever $Q_1,Q_2\in\mathcal{A}$ with $Q_1 \ne Q_2$.
If $\bigcup_{Q\in\mathcal{A}}Q=Q_0$, we say that $\mathcal{A}$ is a \defn{section relative to $Q_0$}, and we say that a section is \defn{complete} if it is a section relative to $\Omega$.
Note that sections are necessarily countable and, for example, each $\mathcal{T}_k$ for $k\in\N\cup\{0\}$ is a section relative to $\Omega$.
The set of sections is equipped with a partial order $\mathcal{A}_1\preccurlyeq\mathcal{A}_2$ if for all $Q_1\in\mathcal{A}_1$ there is a $Q_2\in\mathcal{A}_2$ such that $Q_2\subset Q_1$.
In this situation, we say that $\mathcal{A}_1$ is \emph{refined by} $\mathcal{A}_2$.
This partial order is equipped with a \emph{meet}: that is, given a finite family of sections $\mathcal{A}_1,\ldots,\mathcal{A}_n$, there is a unique section $\mathcal{A}_1\wedge\cdots\wedge\mathcal{A}_n$ which is maximal with respect to the partial order such that
\begin{equation*}
    \mathcal{A}_1\wedge\cdots\wedge\mathcal{A}_n\preccurlyeq\mathcal{A}_i
\end{equation*}
for all $i=1,\ldots,n$.

A metric tree is equipped with a natural family of sections relative to $\Omega$ which respect the geometry of the metric $d$.
We define
\begin{equation}\label{e:geom-section}
    \mathcal{T}(r)=\{Q\in\mathcal{T}:\rho(Q)\leq r<\rho(\widehat{Q})\}
\end{equation}
where, abusing notation, we write $\rho(\widehat{\Omega})=\infty$.
Property \cref{im:rho-sep} above ensures that this is indeed a section and property \cref{im:rho-finite} ensures that $\mathcal{T}_k\preccurlyeq\mathcal{T}(r)$ for all $k$ sufficiently large.

\subsection{Regularity of non-autonomous self-similar sets}\label{sec:non-autonomous}
We now recall some of the results on non-autonomous self-similar sets from \cite{zbmath:08103993}.
We require more than just the dimension of non-autonomous self-similar sets (as already follows, for instance, from \cite{zbl:1364.28011}): we will require additional regularity properties of the similarity dimensions of the finite components of the non-autonomous IFS.

For each $n\in\N$, let $\mathcal{J}_n$ be a finite non-empty index set, and let $\Phi_n=\{S_{n,j}\}_{j\in\mathcal{J}_n}$ be a family of similarity maps $S_{n,j}\colon\R^d\to\R^d$ of the form
\begin{equation*}
    S_{n,j}(\bm{x}) = r_{n,j} O_{n,j}\bm{x} + \bm{d}_{n,j}
\end{equation*}
where $r_{n,j}\in(0,1)$ and $O_{n,j}$ is an orthogonal matrix.
To avoid degenerate situations, we assume that
\begin{equation} \label{e:rn-lim}
    \lim_{n\to\infty}\sup\{r_{1,j_1}\cdots r_{n,j_n}:j_i\in\mathcal{J}_i\text{ for each }i=1,\ldots,n\}=0
\end{equation}
and that there is a compact set $X\subset\R^d$ which is \emph{invariant} such that $S_{n,j}(X)\subset X$ for all $n\in\N$ and $j\in\mathcal{J}_n$.
Associated with the sequence $(\Phi_n)_{n=1}^\infty$, there is a non-empty and compact \defn{limit set}
\begin{equation*}
    K=\bigcap_{n=1}^\infty\bigcup_{(j_1,\ldots,j_n)\in\mathcal{J}_1\times\cdots\times\mathcal{J}_n} S_{1,j_1}\circ\cdots\circ S_{n,j_n}(X).
\end{equation*}
Under these assumptions, the sequence $(\Phi_n)_{n=1}^\infty$ is called a \emph{non-autonomous iterated function system (IFS)} and the limit set $K$ is called the \emph{non-autonomous self-similar set}.

We can associated a natural metric tree to every non-autonomous self-similar set.
Let $\mathcal{T}=\bigcup_{n=0}^\infty\mathcal{T}_n$ denote the set of all cylinder sets in the infinite product space $\Delta=\prod_{n=1}^\infty\mathcal{J}_n$, where
\begin{equation*}
    \mathcal{T}_n = \left\{[j_1,\ldots,j_n] = \{j_1\}\times\cdots\times\{j_n\}\times\prod_{k=n+1}^\infty \mathcal{J}_k:j_i \in\mathcal{J}_i\quad\text{for}\quad i=1,\ldots,n\right\}.
\end{equation*}
Note that the unique cylinder in $\mathcal{T}_0$ is the set $\Delta$.
Given a cylinder $Q=[j_1,\ldots,j_n]\in\mathcal{T}$, we write $\rho(Q) = r_{1,j_1}\cdots r_{n,j_n}$.
The triplet $(\Delta,\{\mathcal{T}_k\}_{k=0}^\infty,\rho)$ is then a metric tree.
We define $\pi\colon\Delta\to\R^d$ by the relation
\begin{equation*}
    \{\pi\bigl((i_n)_{n=1}^\infty\bigr)\} = \bigcap_{n=1}^\infty S_{1,i_1}\circ\cdots\circ S_{n,i_n}(X).
\end{equation*}
This function is well-defined by \cref{e:rn-lim}, and it is easy to see that $\pi$ is Lipschitz so that $\pi(\Delta)=K$.

We say that a non-autonomous IFS $(\Phi_n)_{n=1}^\infty$ satisfies the \defn{open set condition} if the invariant compact set $X$ can be chosen to have non-empty interior $U=X^\circ$ so that for each $n\in\N$ and $j\neq j'\in\mathcal{J}_n$, we have $S_{n,j}(U)\cap S_{n,j'}(U)=\varnothing$.
Furthermore, we say that the IFS has \emph{uniformly bounded contraction ratios} if there are constants $0<r_{\min}\leq r_{\max}<1$ so that $r_{\min}\leq r_{n,j}\leq r_{\max}$ for all $n\in\N$ and $j\in\mathcal{J}_n$.

Our first lemma follows from \cite[Theorem~2.9 and Proposition~2.5]{zbmath:08103993}.
\begin{lemma}\label{c:symb-Assouad}
    Let $(\Phi_n)_{n=1}^\infty$ be a non-autonomous self-similar IFS satisfying the open set condition and with uniformly bounded contraction ratios.
    Denote the associated non-autonomous self-similar set by $K$ and the metric tree by $\Delta$.
    Then $\dimA K=\dimA\Delta$.
\end{lemma}
For each non-autonomous IFS $(\Phi_n)_{n=1}^\infty$ let $\theta(n,m)$ be the similarity dimension of the IFS
$\Phi_{n+1}\circ\cdots\circ\Phi_{n+m}=\{f_1\circ\cdots\circ f_m:f_i\in\Phi_{n+i}\}$
defined by
\begin{equation*}
  \sum_{j_1\in\mathcal{J}_{n}}\cdots\sum_{j_m\in\mathcal{J}_{n+m-1}}\prod_{k=0}^{m-1} r_{n+k,j_{n+k}}^{\theta(n,m)}=1.
\end{equation*}
The following theorem follows by combining \cite[Theorem~2.9, Proposition~3.5, Lemma~4.2, and Theorem~4.3]{zbmath:08103993}.

\begin{theorem}\label{t:non-auto-Assouad}
    Let $(\Phi_n)_{n=1}^\infty$ be a non-autonomous IFS satisfying the open set condition and with uniformly bounded contraction ratios.
    Denote the associated non-autonomous self-similar set by $K$.
    Then for every $n,m,k\in\N$ and $n\leq n'\leq n+k$,
    \begin{enumerate}[nl,r]
        \item\label{i:submax} $\theta(n,m+k) \leq \max\{\theta(n,m),\theta(n+m, k)\}$, and
        \item\label{i:cont} $|\theta(n,m+k) - \theta(n',m)|\lesssim \frac{k}{m}$.
    \end{enumerate}
    Moreover,
    \begin{equation}\label{e:sub-lim}
        \dimA K=\lim_{m\to\infty}\sup_{n\in\N}\theta(n,m)=\lim_{\substack{m\to\infty\\m\in \kappa\N}}\sup_{n\in \kappa\N}\theta(n,m)
    \end{equation}
    for all $\kappa\in\N$, where $\kappa\N=\{\kappa n:n\in\N\}$.
\end{theorem}
Note in \cite[Lemma~4.2]{zbmath:08103993} the inequality corresponding to \cref{i:cont} is one-sided, but in our case the reverse inequality follows from a short convexity argument using $r_{\min}\leq r_{n,j}\leq r_{\max}$.
\begin{remark}
    The observant reader will observe that we also allow index sets with $\#\mathcal{J}_n=1$, whereas \cite{zbmath:08103993} requires that $\#\mathcal{J}_n \geq 2$.
    The assumption $\#\mathcal{J}_n\geq 2$ is only used in the proof of \cite[Lemma~2.7]{zbmath:08103993}, to deduce that products of contraction ratios over $m$ consecutive levels converge to zero uniformly in the starting level as $m\to\infty$.
    In our case, this follows directly from the upper bound $r_{n,j}\leq r_{\max}<1$.
    When $\#\mathcal{J}_n=1$ for some $n$, we interpret cylinders as words rather than as subsets of $\Delta$, so that $\rho$ remains well defined even though consecutive sections may repeat.
    The identity $\diam(Q)=\rho(Q)$ from \cref{ss:metric} weakens to the inequality $\diam(Q)\leq\rho(Q)$, but only the weaker inequality is required.
\end{remark}

\section{Tangent structure and dimension of Gatzouras--Lalley carpets}\label{s:carpet-tangents}
In this section, we introduce the definitions of Gatzouras--Lalley and Barański carpets and prove our main results on tangents and pointwise Assouad dimension of Gatzouras--Lalley carpets.
We defer the results concerning Barański carpets to \cref{sec:baranski}.

We begin by recalling the basic definitions and dimension theory for these carpets in \cref{ss:carpet-def} and \cref{ss:approx-square}.
Next, under additional regularity assumptions, we obtain pointwise bounds on dimension.
The lower bounds on pointwise Assouad dimension in \cref{ss:tangents} are obtained by constructing large tangents corresponding to products of the projection with weak tangents of the non-autonomous self-similar sets corresponding to fibres.
This uses the infrastructure established in \cref{sec:non-autonomous}; note that the lower bounds do not require regularity.
Pointwise upper bounds are then proven in \cref{ss:ub} using a fine covering arguments concerning approximate squares.
Here, the regularity assumptions are essential since we obtain bounds in terms of the symbolic coding of the point.
Finally, we isolate large uniformly regular subsystems with prescribed symbolic dimensions, completing the proof of \cref{it:gl-results}.

\subsection{Gatzouras--Lalley and Barański carpets}\label{ss:carpet-def}
\subsubsection{Defining the maps}
Fix an index set $\mathcal{I}$ with $\#\mathcal{I}\geq 2$, and for $j=1,2$ fix contraction ratios $(\ctr_{i,j})_{i\in\mathcal{I}}\subset(0,1)$ and translations $(d_{i,j})_{i\in\mathcal{I}}\subset\R$.
We then call the IFS $\{T_i\}_{i\in\mathcal{I}}$ \defn{diagonal} when
\begin{equation*}
    T_i(x_1,x_2)=(\ctr_{i,1} x_1+d_{i,1},\ctr_{i,2} x_2+d_{i,2})\quad\text{for each}\quad i\in\mathcal{I}.
\end{equation*}
Let $\eta_j$ denote the orthogonal projection onto the $j$\textsuperscript{th} coordinate axis, i.e.~$\eta_j(x_1,x_2)=x_j$.
We denote by $\Lambda_j=\{S_{i,j}\}_{i\in\mathcal{I}}$ the \defn{projected systems}, where $\eta_j\circ T_j=S_{i,j}\circ\eta_j$.
We will often write $\eta=\eta_1$ to denote simply the projection onto the first coordinate axis.
Of course, $S_{i,j}(x)=\ctr_{i,j} x+d_{i,j}$ are iterated function systems of similarities.

Let $\mathcal{I}^*=\bigcup_{n=0}^\infty\mathcal{I}^n$, and for $\mtt{i}=(i_1,\ldots,i_n)\in\mathcal{I}^*$ and $j=1,2$, write
\begin{align*}
    T_\mtt{i} &= T_{i_1}\circ\cdots\circ T_{i_n},\\*
    S_{\mtt{i},j} &= S_{i_1,j}\circ\cdots\circ S_{i_n,j}
\end{align*}
and
\begin{align*}
    p_\mtt{i} &= p_{i_1}\cdots p_{i_n},\\*
    \ctr_{\mtt{i},j} &= \ctr_{i_1,j}\cdots \ctr_{i_n,j}.
\end{align*}
For $n\in\N$ and $\gamma\in\Omega \coloneqq \mathcal{I}^{\N}$, we write $\gamma\npre{n}$ to denote the unique prefix of $\gamma$ in $\mathcal{I}^{n}$.

Now $\eta_j$ induces an equivalence relation $\sim_j$ on $\mathcal{I}$ where $i\sim_j i'$ if $S_{i,j}=S_{i',j}$.
Let $\eta_j\colon\mathcal{I}\to\mathcal{I}/\sim_j$ denote the natural projection.
Intuitively, $\eta_j(i)$ is the set of indices which lie in the same column or row as the index $i$.
Then $\eta_j$ extends naturally to a map on $\Omega$ by $\eta_j((i_n)_{n=1}^\infty)=(\eta_j(i_n))_{n=1}^\infty\subset\eta_j(\mathcal{I})^{\N}\cong\eta_j(\mathcal{I}^{\N})$; and similarly extends to a map on $\mathcal{I}^*$.
For notational clarity, we will refer to words in $\mathcal{I}^*$ using upright indices, such as $\mtt{i}$, and words in $\eta_j(\mathcal{I}^*)$ using their underlined variants, such as $\umtt{i}$.
Note that if $\mtt{i}\sim_j\mtt{j}$, then and $S_{\mtt{i},j}=S_{\mtt{j},j}$.
In particularly, we may unambiguously write $S_{\umtt{i},j}$ and $\ctr_{\umtt{i},j}$ for $\umtt{i}\in\eta_j(\mathcal{I}^*)$.

Associated with the IFS $\{T_i\}_{i\in\mathcal{I}}$ is a unique non-empty compact \defn{attractor} $K$, satisfying $K=\bigcup_{i\in\mathcal{I}}T_i(K)$.
Note that the projected IFS $\{S_{i,j}\}_{i\in\mathcal{I}}$ has attractor $K_j=\eta_j(K)$ for $j=1,2$.
Recalling that $\Omega=\mathcal{I}^{\N}$, let $\pi\colon\Omega\to K$ denote the continuous map uniquely defined by
\begin{equation*}
    \bigl\{\pi\bigl((i_n)_{n=1}^\infty\bigr)\bigr\}=\lim_{n\to\infty}S_{i_1}\circ\cdots\circ S_{i_n}(K).
\end{equation*}
Without loss of generality, and for the remainder of this document, we will assume that $K\subset[0,1]^2$.
We can now introduce our two primary classes of self-affine sets.
\begin{definition}
    We say that the carpet is of type \emph{Gatzouras--Lalley} if:
    \begin{enumerate}[nl]
        \item $T_i((0,1)^2)\cap T_j((0,1)^2)=\varnothing$ for all $i\neq j$,
        \item either $S_{i,1}((0,1))=S_{(j,1)}((0,1))$ or $S_{i,1}((0,1))\cap S_{j,1}((0,1))=\varnothing$ for all $i,j$, and
        \item $\ctr_{i,1}>\ctr_{i,2}$ for all $i\in\mathcal{I}$;
    \end{enumerate}
    and type \emph{Barański} if:
    \begin{enumerate}[nl]
        \item $T_i((0,1)^2)\cap T_j((0,1)^2)=\varnothing$ for all $i\neq j$, and
        \item either $S_{i,\ell}((0,1))=S_{(j,\ell)}((0,1))$ or $S_{i,\ell}((0,1))\cap S_{j,\ell}((0,1))=\varnothing$ for all $i,j$ and $\ell=1,2$.
    \end{enumerate}
\end{definition}
Moreover, we say that an IFS $\{f_i\}_{i\in\mathcal{I}}$ with attractor $K$ satisfies the \defn{strong separation condition} (or SSC for short) if $f_i(K)\cap f_j(K)=\varnothing$ for all $i\neq j\in\mathcal{I}$.

\subsubsection{Dimensions of Gatzouras--Lalley carpets}\label{sss:dims}
To conclude this section, we recall some standard results on the dimensions of Gatzouras--Lalley carpets.
We defer the corresponding results for Barański carpets to \cref{sec:baranski-dim}.

Before we do this, we first recall the notion of the lower dimension, which is in some sense dual to the definition of Assouad definition.
Let $K\subset\R^d$ be compact.
Then the \defn{lower dimension} of $K$ is given by
\begin{align*}
    \dimL K=\sup\Bigl\{s:\exists C>0\,&\forall 0<r\leq R<1\,\forall x\in K\\
                                      &N_r(B(x,R)\cap K)\geq C\Bigl(\frac{R}{r}\Bigr)^s\Bigr\}.
\end{align*}
In order to state our results on the Hausdorff dimensions, we must also introduce some notation for Bernoulli measures.
Let $\mathcal{P}$ denote the collection of probability vectors on $\mathcal{I}$, i.e.
\begin{equation*}
    \mathcal{P}=\mathcal{P}(\mathcal{I})\coloneqq\Bigl\{(p_i)_{i\in\mathcal{I}}:p_i\geq 0\text{ for all }i\text{ and }\sum_{i\in\mathcal{I}}p_i=1\Bigr\}.
\end{equation*}
Equip $\mathcal{P}$ with the topology inherited from $\R^{\mathcal{I}}$.
Given $\bm{p}\in\mathcal{P}$, considering $\bm{p}$ as a probability measure on $\mathcal{I}$, we let $\bm{p}^{\N}$ denote the infinite product measure supported on $\Omega$.
We let $\mu_{\bm{p}}=\pi_*\bm{p}^{\N}$ denote the corresponding invariant measure on $K$, where $\pi_*$ denotes the pushforward map.
Note that the projections $\eta_j$ also induce natural maps $\eta_j\colon\mathcal{P}(\mathcal{I})\to\mathcal{P}(\eta_j(\mathcal{I}))$ by $\eta_j(\bm{p})_{\underline{\ell}}=\sum_{i\in\eta_j^{-1}(\underline{\ell})}\bm{p}_i$.

Given a probability vector $\bm{p}\in\mathcal{P}$, we write
\begin{equation*}
    H(\bm{p})=-\sum_{i\in\mathcal{I}}p_i\log p_i\qquad\text{and}\qquad\chi_j(\bm{p})=-\sum_{i\in\mathcal{I}}p_i\log \ctr_{i,j}.
\end{equation*}
We now recall the main results of \cite{zbl:0757.28011}---stated below in \cref{im:gl-hdim} and \cref{im:gl-bdim}---as well as the result of \cite{zbl:1278.37032}---stated below in \cref{im:gl-adim}.
We also note that the same proof as given in \cite{zbl:1278.37032} (which is explained more precisely in \cite[Theorem~2.13]{zbl:1305.28021}) gives the analogous result for the lower dimension.
\begin{proposition}[\cite{zbl:0757.28011,zbl:1278.37032}]\label{p:gl-dims}
    Let $K$ be a Gatzouras--Lalley carpet.
    \begin{enumerate}[r]
        \item\label{im:gl-hdim} The Hausdorff dimension of $K$ is given by
            \begin{equation*}
                \dimH K=\sup_{\bm{p}\in\mathcal{P}}s(\bm{p})
            \end{equation*}
            where
            \begin{equation*}
                s(\bm{p})\coloneqq\frac{H(\eta(\bm{p}))}{\chi_1(\bm{p})}+\frac{H(\bm{p})-H(\eta(\bm{p}))}{\chi_2(\bm{p})}.
            \end{equation*}
            Moreover, the supremum is always attained at an interior point of $\mathcal{P}$ (i.e.\ at vector $\bm{w}\in\mathcal{P}$ with $\bm{w}_i>0$ for all $i\in\mathcal{I}$).

        \item\label{im:gl-bdim} The box dimension of $K$ exists and is given by the unique solution to
            \begin{equation*}
                \sum_{i\in\mathcal{I}}\ctr_{i,1}^{\dimB\eta(K)}\ctr_{i,2}^{\dimB K-\dimB\eta(K)}=1\qquad\text{where}\qquad\sum_{\underline{j}\in\eta(\mathcal{I})}\ctr_{\underline{j},1}^{\dimB\eta(K)}=1.
            \end{equation*}

        \item\label{im:gl-adim} The Assouad dimension of $K$ is given by
            \begin{equation*}
                \dimA K=\dimB\eta(K)+\max_{\underline{\ell}\in\eta(\mathcal{I})} t(\underline{\ell})
            \end{equation*}
            where $t(\underline{\ell})$ is defined as the unique solution to the equations
            \begin{equation*}
                \sum_{j\in\eta^{-1}(\underline{\ell})}\ctr_{j,2}^{t(\underline{\ell})}=1.
            \end{equation*}
            Similarly, the lower dimension of $K$ is given by
            \begin{equation*}
                \dimL K=\dimB\eta(K)+\min_{\underline{\ell}\in\eta(\mathcal{I})} t(\underline{\ell}).
            \end{equation*}
    \end{enumerate}
\end{proposition}

\subsubsection{Regular points and interior words}
We conclude this section with the notion of a regular point and an interior word.
Heuristically, a regular point is one which uniformly avoids points in other columns; note that we do not require such a property \emph{within} each columns.
\begin{definition}\label{d:regular}
    We say that a point $x\in K$ is \emph{regular} if for each $r\in(0,1)$, there is an $\mtt{i}\in\mathcal{I}^*$ with $\ctr_{\mtt{i},1}\lesssim r$ such that $B(\eta(x),r)\cap\eta(K)\subset S_{\mtt{i},1}(\eta(K))$.
    Given $\mtt{i}\in\mathcal{I}^*$, we say that $\mtt{i}$ is an \defn{interior word} if $S_{\mtt{i},1}([0,1])\subset(0,1)$.
    We let $\mathcal{B}_n\subset\mathcal{I}^n$ denote the set of interior words of length $n$.
\end{definition}

The following lemma is standard.
Recall that $\Omega=\mathcal{I}^{\N}$ is the symbolic space coding the attractor $K$.
Here, and elsewhere, given an $n\in\N$ and $\mathcal{Y}\subset\mathcal{I}^n$, we embed $\mathcal{Y}^{\N}$ in $\Omega$ in the natural way.
We will abuse notation and interchangeably refer to elements in the subsystem or in the full system.
\begin{lemma}\label{l:gl-reg}
    Let $K$ be a Gatzouras--Lalley carpet.
    \begin{enumerate}[nl,r]
        \item\label{im:vssc} If $\eta(K)$ satisfies the SSC, then each $x\in K$ is regular.
        \item\label{im:all-subword} Suppose $\gamma\in\mathcal{B}_n^{\N}$ for some $n\in\N$.
            Then $\pi(\gamma)$ is regular.
    \end{enumerate}
\end{lemma}
We can now guarantee the existence of large subsystems consisting only of regular points.
This result is essentially \cite[Lemma~4.3]{zbl:1206.28011}.
\begin{proposition}[\cite{zbl:1206.28011}]\label{p:gl-large-subsystem}
    Let $K$ be a Gatzouras--Lalley carpet corresponding to the IFS $\{T_i\}_{i\in\mathcal{I}}$.
    Then for every $\epsilon>0$, there is an $n\in\N$ and a family $\mathcal{J}\subset\mathcal{I}^n$ so that the IFS $\{T_{\mtt{j}}:\mtt{j}\in\mathcal{J}\}$ with attractor $K_\epsilon$ satisfies the following conditions:
    \begin{enumerate}[nl,r]
        \item\label{im:interior} each $\mtt{i}\in\mathcal{J}$ is an interior word,
        \item\label{im:dim} $\dimH K_\epsilon\geq\dimH K-\epsilon$,
        \item\label{im:dim-proj} $\dimB\eta(K_\epsilon)\geq \dimB\eta(K)-\epsilon$, and
        \item\label{im:uniform} there are $0<\rho_2<\rho_1<1$ so that $\ctr_{\mtt{i},1}=\rho_1$ and $\ctr_{\mtt{i},2}=\rho_2$ for all $\mtt{i}\in\mathcal{I}$ and each column has the same number of maps.
    \end{enumerate}
    In particular, each $x\in K_\epsilon$ is a regular point with respect to the IFS $\{T_i\}_{i\in\mathcal{I}}$ and $\dimA K_\epsilon =\dimH K_\epsilon=\dimL K_\epsilon$.
\end{proposition}
\begin{proof}
    First, if $K$ is contained in a vertical line, then $K$ is the attractor of a self-similar IFS in $\R$ and the result is substantially easier.
    Now applying \cite[Lemma~4.3]{zbl:1206.28011}, there exists a family $\mathcal{J}_0\subset\mathcal{I}^{n_0}$ with attractor $K_0$ satisfying conditions \cref{im:dim}, \cref{im:dim-proj}, and \cref{im:uniform}.
    By condition \cref{im:uniform}, there is a $t\in\R$ so that $t(\mtt{i})=t$ for all $\mtt{i}\in\mathcal{J}_0$.
    Therefore
    \begin{equation*}
        \dimH K_0 = \dimB\eta(K)+t
    \end{equation*}
    and since $K$ is not contained in a vertical line, we may assume that $\dimB\eta(K_0)>0$.

    Since $\eta(K_0)$ is the attractor of a self-similar IFS, iterating $\mathcal{J}_0$ if necessary and removing the maps in the first and last column, obtain a family $\mathcal{J}\subset\mathcal{J}_0^n$ with corresponding attractor $K_\epsilon$ such that $t(\mtt{j})=t$ for any $\mtt{j}\in\mathcal{J}$, and $\dimB\eta(K_\epsilon)\geq\dimB\eta(K)-\epsilon$.
    Since words which correspond to rectangles that do not lie in the first or last column are necessarily interior words, combining this construction with \cref{l:gl-reg} provides a family $\mathcal{J}$ satisfying the desired properties.
\end{proof}

\subsection{Approximate squares and symbolic slices}\label{ss:approx-square}
A common technique when studying invariant sets for iterated function systems on some index set $\mathcal{I}$ is to first reduce the problem to a symbolic problem on the coding space $\mathcal{I}^*$.
However, the main technical complexity in understanding the dimension theory Gatzouras--Lalley carpets, and more generally self-affine sets, is that the cylinder sets $T_{\mtt{i}}(K)$ are often exponentially distorted rectangles.
As a result, we will keep track of two symbolic systems simultaneously, which together will capture the geometry of the set $K$.

Fix a Gatzouras--Lalley IFS $\Lambda=\{T_i\}_{i\in\mathcal{I}}$.
We first introduce some notation for handling cylinders.
We then associate with the IFS $\Lambda$, and the related defining data that we introduced in \cref{ss:carpet-def}, two important metric trees: first, the metric tree of \defn{approximate squares}, and second the metric tree of \defn{symbolic slices}.

First, recall that $\Omega=\mathcal{I}^{\N}$ is the space of infinite sequences on $\mathcal{I}$.
Given $k\in\N\cup\{0\}$ and a word $\mtt{i}\in\mathcal{I}^k$, we define the \defn{cylinder} corresponding to $\mtt{i}$ by
\begin{equation*}
    [\mtt{i}]=\{\gamma\in\Omega:\gamma\npre{k}=\mtt{i}\}.
\end{equation*}
The family of cylinders $\{[\mtt{i}]:\mtt{i}\in\mathcal{I}^k\}_{k=0}^\infty$ defines a tree: we will often abuse notation and simply refer to $\{\mathcal{I}^k\}_{k=0}^\infty$ as a tree.
We will associate with this tree a variety of metrics, such as those induced by the maps $\mtt{i}\mapsto\ctr_{\mtt{i},j}$ for $j=1,2$.
We will also use the same notation for the projected words $\{\eta(\mathcal{I}^k)\}_{k=0}^\infty$.

Next, we define the \defn{metric tree of approximate squares}.
Before we do this, we introduce the notion of a \defn{pseudo-cylinder}.
\begin{figure}[t]
    \centering
    \begin{tikzpicture}[scale=2.9]
    \draw[fill=blue!10, opacity=0.5] (0,0) rectangle (5/2, 1/3+1/4);
    \draw[fill=red!10, opacity=0.5] (3*5/4,0) rectangle (5,1);

    \begin{scope}[thick]
        \draw (0.0, 0.0) -- (2.5, 0.0) -- (2.5, 0.25) -- (0.0, 0.25) -- cycle;
        \draw (0.0, 0.3333333333333333) -- (2.5, 0.3333333333333333) -- (2.5, 0.5833333333333333) -- (0.0, 0.5833333333333333) -- cycle;
        \draw (2.5, 0.0) -- (5.0, 0.0) -- (5.0, 0.2) -- (2.5, 0.2) -- cycle;
        \draw (2.5, 0.6666666666666666) -- (5.0, 0.6666666666666666) -- (5.0, 1.0) -- (2.5, 1.0) -- cycle;
    \end{scope}

    \begin{scope}
        \draw (0.0, 0.0) -- (1.25, 0.0) -- (1.25, 0.0625) -- (0.0, 0.0625) -- cycle;
        \draw (0.0, 0.3333333333333333) -- (1.25, 0.3333333333333333) -- (1.25, 0.3958333333333333) -- (0.0, 0.3958333333333333) -- cycle;
        \draw (2.5, 0.0) -- (3.75, 0.0) -- (3.75, 0.05) -- (2.5, 0.05) -- cycle;
        \draw (2.5, 0.6666666666666666) -- (3.75, 0.6666666666666666) -- (3.75, 0.75) -- (2.5, 0.75) -- cycle;
        \draw (0.0, 0.08333333333333333) -- (1.25, 0.08333333333333333) -- (1.25, 0.14583333333333331) -- (0.0, 0.14583333333333331) -- cycle;
        \draw (0.0, 0.41666666666666663) -- (1.25, 0.41666666666666663) -- (1.25, 0.47916666666666663) -- (0.0, 0.47916666666666663) -- cycle;
        \draw (2.5, 0.06666666666666667) -- (3.75, 0.06666666666666667) -- (3.75, 0.11666666666666665) -- (2.5, 0.11666666666666665) -- cycle;
        \draw (2.5, 0.7777777777777777) -- (3.75, 0.7777777777777777) -- (3.75, 0.861111111111111) -- (2.5, 0.861111111111111) -- cycle;
        \draw (1.25, 0.0) -- (2.5, 0.0) -- (2.5, 0.05) -- (1.25, 0.05) -- cycle;
        \draw (1.25, 0.3333333333333333) -- (2.5, 0.3333333333333333) -- (2.5, 0.3833333333333333) -- (1.25, 0.3833333333333333) -- cycle;
        \draw (3.75, 0.0) -- (5.0, 0.0) -- (5.0, 0.04) -- (3.75, 0.04) -- cycle;
        \draw (3.75, 0.6666666666666666) -- (5.0, 0.6666666666666666) -- (5.0, 0.7333333333333333) -- (3.75, 0.7333333333333333) -- cycle;
        \draw (1.25, 0.16666666666666666) -- (2.5, 0.16666666666666666) -- (2.5, 0.25) -- (1.25, 0.25) -- cycle;
        \draw (1.25, 0.5) -- (2.5, 0.5) -- (2.5, 0.5833333333333333) -- (1.25, 0.5833333333333333) -- cycle;
        \draw (3.75, 0.13333333333333333) -- (5.0, 0.13333333333333333) -- (5.0, 0.2) -- (3.75, 0.2) -- cycle;
        \draw (3.75, 0.8888888888888888) -- (5.0, 0.8888888888888888) -- (5.0, 1.0) -- (3.75, 1.0) -- cycle;
    \end{scope}
\end{tikzpicture}
    \caption{Two iterations of a Gatzouras--Lalley IFS within a cylinder, with a wide pseudo-cylinder in highlighted in blue and a tall pseudo-cylinder in red.}
    \label{f:pseudo-cylinder}
\end{figure}
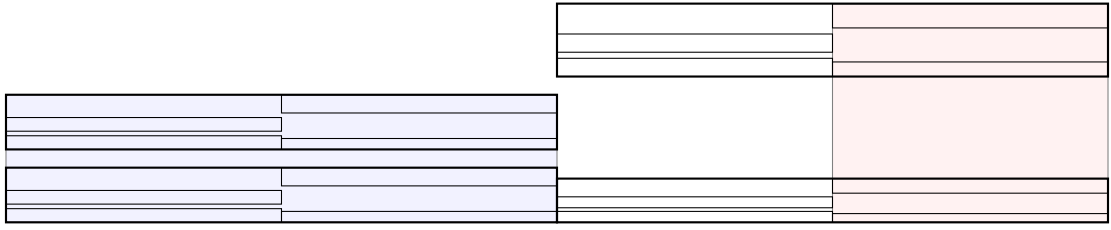
Suppose $\mtt{i}\in\mathcal{I}^k$ and $\umtt{j}\in\eta(\mathcal{I}^\ell)$.
We then write
\begin{equation*}
    P(\mtt{i},\umtt{j})=\{\gamma=(i_n)_{n=1}^\infty\in\Omega:(i_1,\ldots,i_k)=\mtt{i}\text{ and }\eta(i_{k+1},\ldots,i_{k+\ell})=\umtt{j}\}.
\end{equation*}
Note that map $(\mtt{i},\umtt{j})\mapsto P(\mtt{i},\umtt{j})$ is injective.
Another equivalent way to understand the pseudo-cylinder $P(\mtt{i},\umtt{j})$ is as a finite union of cylinders inside the cylinder $[\mtt{i}]$, all of which lie inside the same column; that is,
\begin{equation}\label{e:pseudo-cylinder-rep}
    P(\mtt{i},\umtt{j})=\bigcup_{\mtt{k}\in\eta^{-1}(\umtt{j})}[\mtt{i}\mtt{k}].
\end{equation}
We refer the reader to \cref{f:pseudo-cylinder} for a depiction of the definition of a pseudo-cylinder.

Now given an infinite word $\gamma\in\Omega$, let $L_k(\gamma)$ be the minimal integer so that
\begin{equation*}
    \ctr_{\gamma_1,1}\cdots\ctr_{\gamma_{L_k(\gamma)},1}<\ctr_{\gamma_1,2}\cdots\ctr_{\gamma_k,2}.
\end{equation*}
In other words, $L_k(\gamma)$ is chosen so that the level $L_k(\gamma)$ rectangle has approximately the same width as the height of the level $k$ rectangle.
Write $\gamma\npre{L_k(\gamma)}=\mtt{i}\mtt{j}$ where $\mtt{i}\in\mathcal{I}^k$.
We then define the \defn{approximate square} $Q_k(\gamma)\subset\Omega$ by
\begin{equation*}
    Q_k(\gamma)=P(\mtt{i},\eta(\mtt{j})).
\end{equation*}
While different $\gamma$ may define the same approximate square, the choice of $\mtt{i}$ and $\eta(\mtt{j})$ are unique.
For fixed $\mtt{i}$, let $\mathcal{U}(\mtt{i})\subset\eta(\mathcal{I}^*)$ denote the set of $\umtt{j}$ so that $P(\mtt{i},\umtt{j})$ is an approximate square.
Of course, $Q_{k+1}(\gamma)\subset Q_k(\gamma)$ and moreover for any $\gamma,\gamma'\in\Omega$, either $Q_{k}(\gamma)=Q_{k}(\gamma')$ or $Q_k(\gamma)\cap Q_k(\gamma')=\varnothing$.
In particular, $\mathcal{U}(\mtt{i})$ is a complete section and the approximate squares $P(\mtt{i},\umtt{j})$ are disjoint in symbolic space for fixed $\mtt{i}$.

We say that a pseudo-cylinder $P(\mtt{i},\umtt{j})$ is \emph{wide} if $\umtt{j}\preccurlyeq\umtt{k}$ for some $\umtt{k}\in\mathcal{U}(\mtt{i})$; in other words, $P(\mtt{i},\umtt{j})$ contains approximate squares of the form $P(\mtt{i},\umtt{k})$.
Otherwise, we say that $P(\mtt{i},\umtt{j})$ is \emph{tall}.
In other words, one can think of the wide pseudo-cylinders as ``interpolating'' between the cylinder $P(\mtt{i},\varnothing)=[\mtt{i}]$ and the approximate square $P(\mtt{i},\umtt{j})=Q_n(\gamma)$.

Denote the set of all approximate squares by
\begin{equation*}
    \mathcal{S}_k=\{Q_k(\gamma):\gamma\in\Omega\}\qquad\text{and}\qquad\mathcal{S}=\bigcup_{k=0}^\infty\mathcal{S}_k.
\end{equation*}
As discussed above, every approximate square is uniquely associated with a pair $(\mtt{i},\umtt{j})$, so we may therefore define a metric induced by $\rho(Q)=\ctr_{\mtt{i},2}$, which makes the collection of approximate squares into a metric tree.

To conclude this section, we define the \defn{metric tree of symbolic slices}.
Suppose we fix a word $\gamma\in\Omega$.
The word $\gamma=(i_n)_{n=1}^\infty$ defines for each $n\in\N$ a self-similar IFS $\Phi_n=\{S_{j,2}:j\in\eta^{-1}(\eta(i_n))\}$.
This IFS is precisely the IFS corresponding to the column containing the index $i_n$.
Note that there are only finitely many possible choices for the $\Phi_n$, so the sequence $(\Phi_n)_{n=1}^\infty$ has as an attractor a non-autonomous self-similar set $K_{\eta(\gamma)}$ and corresponding metric tree $\Omega(\eta(\gamma))$, as defined in \cref{sec:non-autonomous}.
This non-autonomous IFS has uniformly bounded contractions and satisfies the OSC with respect to the open interval $(0,1)$.
For notational simplicity, we denote the cylinder sets which compose this metric tree as
\begin{equation*}
    \mathcal{F}_{\eta(\gamma),n}=\{[j_1,\ldots,j_n]:(j_1,\ldots,j_n)\in \Phi_1\times\cdots\times\Phi_n\}\quad\text{and}\quad\mathcal{F}_{\eta(\gamma)}=\bigcup_{n=0}^\infty\mathcal{F}_{\eta(\gamma),n}.
\end{equation*}
We call $K_{\eta(\gamma)}$ the \defn{symbolic slice} associated with the word $\gamma$.
If the projected IFS $\{S_{\underline{i},1}\}_{\underline{i}\in\eta(\mathcal{I})}$ satisfies the SSC, then if $x=\eta(\pi(\gamma))$,
\begin{equation*}
    \{x\}\times K_{\eta(\gamma)}=\eta^{-1}(x)\cap K
\end{equation*}
is precisely the vertical slice of $K$ containing $x$.
In general, $K_{\eta(\gamma)}$ is always contained inside a vertical slice of $K$.
The symbolic fibre $K_{\eta(\gamma)}$ (and its associated Assouad dimension) was introduced and studied in \cite[§1.2]{zbl:1561.28063} in the more general setting of overlapping diagonal carpets.

\subsection{Tangents of Gatzouras--Lalley carpets}\label{ss:tangents}
It turns out that the pointwise Assouad dimension at $x=\pi(\gamma)$ is closely related to the Assouad dimension of the symbolic fibre $K_{\eta(\gamma)}$.
In this section, we make this notion precise, and moreover use it to construct large tangents for Gatzouras--Lalley carpets.

In our main result in this section, we prove that approximate squares containing a fixed word $\gamma\in\Omega$ converge in Hausdorff distance to product sets of weak tangents of $K_{\eta(\gamma)}$ with the projection $\eta(K)$, up to some finite distortion and contributions from adjacent approximate squares.
First, we define
\begin{equation*}
    \Phi_{k,\gamma}(x,y)=\bigl(S_{\gamma\npre{L_k(\gamma)},1}^{-1}(x),S_{\gamma\npre{k},2}^{-1}(y)\bigr).
\end{equation*}
By choice of $L_k(\gamma)$, the maps $\Phi_{k,\gamma}$ are (up to some constant-size distortion) homotheties.
One can think of $\Phi_{k,\gamma}$ as mapping the approximate square $\pi(Q_k(\gamma))$ to the unit square $[0,1]^2$.
\begin{proposition}\label{l:square-structure}
    Let $K$ be a Gatzouras--Lalley carpet and let $\gamma\in\Omega$ be arbitrary.
    Suppose $(\mtt{i}_n)_{n=1}^\infty$ is any sequence such that $\eta(\mtt{i}_n)=\eta(\gamma\npre{n})$.
    Then
    \begin{equation}\label{e:approx-sq-lim}
        p_{\mathcal{H}}\left(\eta(K)\times (S_{\mtt{i}_n,2}^{-1}(K_{\eta(\gamma)})\cap[0,1]);\Phi_{n,\gamma}(K)\cap[0,1]^2\right)\lesssim\kappa^n
    \end{equation}
    where $\kappa=\max\bigl\{\frac{\ctr_{i,2}}{\ctr_{i,1}}:i\in\mathcal{I}\bigr\}\in(0,1)$.
    Moreover, suppose $\gamma$ is regular.
    Then for any $\gamma\in \Omega$ and $F\in\Tan(K,\pi(\gamma))$, there is an $E\in\Tan(K_{\eta(\gamma)})$ and a similarity map $h$ so that $h(F)\subset\eta(K)\times E$.
\end{proposition}
\begin{proof}
    We first prove that
    \begin{equation*}
        d_{\mathcal{H}}\left(\eta(K)\times (S_{\mtt{i}_n,2}^{-1}(K_{\eta(\gamma)})\cap[0,1]),\Phi_{n,\gamma}(\pi(Q_n(\gamma)))\right)\lesssim\kappa^n
    \end{equation*}
    Fix $n\in\N$ and write $k=L_n(\gamma)$.
    Let $Q_n(\gamma)=P(\gamma\npre{n},\umtt{j})$ and enumerate $\eta^{-1}(\umtt{j})=\{\mtt{j}_1,\ldots,\mtt{j}_n\}\subset\mathcal{I}^{k-n}$.
    Observe that $\eta(T_{\mtt{j}_i}(K))=S_{\mtt{j}_i,1}(K)$ does not depend on the choice of $i=1,\ldots,m$.
    Now $\Phi_{n,\gamma}(T_{\gamma\npre{n}\mtt{j}_i}(K))$ is contained in the rectangle $\eta(K)\times S_{\mtt{j}_i,2}(K)$.
    Moreover, the rectangle $\eta(K)\times S_{\mtt{j}_i,2}(K)$ has height $\lesssim\kappa^n$.
    Therefore
    \begin{equation}\label{e:proj-close}
        d_{\mathcal{H}}\left(\eta(K)\times \bigcup_{i=1}^m S_{\mtt{j}_i,2}([0,1]),\Phi_{n,\gamma}(Q_n(\gamma))\right)\lesssim\kappa^n.
    \end{equation}
    But approximating the set $S_{\mtt{i}_n,2}([0,1])\cap K_{\eta(\gamma)}$ at level $n$ with cylinders at level $k=L_n(\gamma)$, using the fact that $\eta(\mtt{i}_n)=\eta(\gamma\npre{n})$,
    \begin{equation}\label{e:fibre-close}
        d_{\mathcal{H}}\left(S_{\mtt{i}_n,2}^{-1}(K_{\eta(\gamma)})\cap [0,1],\bigcup_{i=1}^m S_{\mtt{j}_i,2}([0,1])\right)\lesssim\kappa^n.
    \end{equation}
    Combining \cref{e:proj-close} and \cref{e:fibre-close} gives the claim.
    In particular, noting that $Q_n(\gamma)\subset K$ and $\Phi_{n,\gamma}(Q_n(\gamma))\subset[0,1]^2$ gives \cref{e:approx-sq-lim}.

    Now suppose in addition that $x=\pi(\gamma)$ is regular and let $r>0$ be arbitrary.
    Since $x$ is regular, there is an $n\in\N$ with $r\leq \ctr_{\gamma\npre{n},1}\lesssim r$ such that
    \begin{equation*}
        B(x,r)\cap K\subset\bigcup_{j=1}^\ell T_{\mtt{i}_j}(K)
    \end{equation*}
    where
    \begin{equation*}
        \{\mtt{i}_1,\ldots,\mtt{i}_\ell\}=\{\mtt{i}\in\mathcal{I}^n:\eta(\mtt{i})=\eta(\gamma\npre{n})\text{ and }T_{\mtt{i}}(K)\cap B(x,r)\neq\varnothing\}.
    \end{equation*}
    Now exactly as before, each rectangle $T_{\mtt{i}_j}(K)$ has width $\approx r$ and height $\lesssim r\kappa^n$.
    Therefore identifying $x\in K$ with the analogous point $x\in K_{\eta(\gamma)}$, there is a similarity map $h_r$ with contraction ratio in some interval $[1,c]$ for a fixed $c$ depending only on the IFS so that
    \begin{equation*}
        p_{\mathcal{H}}\left(r^{-1}(K-x)\cap B(0,1); h_r(\eta(K))\times r^{-1}(K_{\eta(\gamma)}-x)\right)\lesssim\kappa^n.
    \end{equation*}
    Now suppose $F\in\Tan(K,x)$ so that $F=\lim_{n\to\infty} r_n^{-1}(K-x)\cap B(0,1)$.
    Passing to a subsequence, we may assume that the $h_{r_n}$ have contraction ratios converging to some $r_0\geq 1$.
    Thus passing again to a subsequence, let $F_0=\lim_{n\to\infty}(r_0 r_n)^{-1}(K-x)\cap B(0,1)$.
    Since $r_0\geq 1$, we have $F\subset F_0$.
    Passing again to a subsequence, let
    \begin{equation*}
        \lim_{n\to\infty}(r_0 r_n)^{-1}(K_{\eta(\gamma)}-x)\cap B(0,1)=E\in\Tan(K_{\eta(\gamma)}).
    \end{equation*}
    Thus $r_0^{-1}F\subset F_0\subset\eta(K)\times E$, as claimed.
\end{proof}
To conclude this section, we establish our general result which guarantees the existence of product-like tangents for arbitrary points in Gatzouras--Lalley carpets.
\begin{proposition}\label{p:gl-pointwise-tangents}
    Let $K$ be a Gatzouras--Lalley carpet.
    Then for each $x \in K$, there is an $F\in\Tan(K,x)$ so that
    \begin{equation*}
        \mathcal{H}^{\dimH\eta(K)+\dimA K_{\eta(\gamma)}}(F)\gtrsim 1,
    \end{equation*}
    where $\gamma \in \Omega$ is such that $\pi(\gamma)=x$.
    In particular,
    \begin{equation*}
        \dimA(K,x)\geq\max\{\dimH\eta(K)+\dimA K_{\eta(\gamma)},\dimB K\}.
    \end{equation*}
\end{proposition}
\begin{proof}
    We will construct the set $F$ essentially as a product $\eta(K)\times E$ where $E$ is a \emph{weak} tangent of $K_{\eta(\gamma)}$.
    First, recall from \cref{c:symb-Assouad} that $\dimA K_{\eta(\gamma)}=\dimA\Omega(\eta(\gamma))$.
    Thus by \cref{c:full-tan}, there is a sequence $(n_k)_{k=1}^\infty$ diverging to infinity and words $\mtt{i}_k\in\mathcal{I}^{n_k}$ with $\eta(\mtt{i}_k)=\gamma\npre{n_k}$ such that
    \begin{equation*}
        E\coloneqq\lim_{k\to\infty}S_{\mtt{i}_k,2}^{-1}(K_{\eta(\gamma)})\cap[0,1]
    \end{equation*}
    has $\mathcal{H}^{\dimA K_{\eta(\gamma)}}(E)\gtrsim 1$.

    Thus by \cref{l:square-structure} applied along the sequence $(\mtt{i}_k)_{k=1}^\infty$, since the images $\Phi_{n,\gamma}^{-1}([0,1]^2)$ are rectangles with bounded eccentricity containing $\pi(\gamma)$, there is a tangent $F\in\Tan(K,x)$ containing an image of $\eta(K)\times E$ under a bi-Lipschitz map with constants depending only on $K$.
    But $\eta(K)$ is Ahlfors--David regular so that
    \begin{equation*}
        \mathcal{H}^{\dimH\eta(K)+\dimA K_{\eta(\gamma)}}(F)\geq\mathcal{H}^{\dimH\eta(K)+\dimA K_{\eta(\gamma)}}\bigl(\eta(K)\times E\bigr)\gtrsim 1
    \end{equation*}
    as claimed.
    The result concerning $\dimA(K,x)$ then follows by \cref{p:pointwise-tangent} and \cref{p:tan-lower}.
\end{proof}

\subsection{Upper bounds for the pointwise Assouad dimension}\label{ss:ub}
We now prove our main upper bound for the pointwise Assouad dimension of Gatzouras--Lalley carpets.
As a result of the local inhomogeneity of Gatzouras--Lalley carpets, obtaining good upper bounds requires some care.
We will prove a sequence of lemmas which, morally, provide optimal covers for a variety of symbolic objects: these covers will then be combined to obtain our general upper bound for the pointwise Assouad dimension.

We first show that, as a result of the vertical alignment of their component cylinders, pseudo cylinders can essentially be covered by their projection.
Recall that $\mathcal{S}$ denotes the set of all approximate squares.
Then if $P(\mtt{i},\umtt{j})$ is any wide pseudo-cylinder, we can write it as a union of the approximate squares in the family
\begin{equation*}
    \mathcal{Q}(\mtt{i},\umtt{j})=\{Q\in\mathcal{S}:\;Q=P(\mtt{i},\umtt{k})\text{ for some }\umtt{k}\in\eta(\mathcal{I}^*)\text{ and }Q\subset P(\mtt{i},\umtt{j})\}.
\end{equation*}
Since each $Q=P(\mtt{i},\umtt{k})$ for some $\umtt{k}$, we have $Q\in\mathcal\mathcal{S}(\ctr_{\mtt{i},2})$ so that this family of approximate squares forms a section.
\begin{lemma}\label{l:gl-pseudo-cyl-cover}
    Let $P(\mtt{i},\umtt{j})$ be a wide pseudo-cylinder.
    Then
    \begin{equation*}
        \#\mathcal{Q}(\mtt{i},\umtt{j})\approx\left(\frac{\ctr_{\mtt{i}\umtt{j},1}}{\ctr_{\mtt{i},2}}\right)^{\dimB\eta(K)}.
    \end{equation*}
\end{lemma}
\begin{proof}
    First, enumerate $\mathcal{Q}(\mtt{i},\umtt{j})=\{Q_1,\ldots,Q_m\}$, and for each $i=1,\ldots,m$, there is a unique $\umtt{k}_i$ so that $Q_i=P(\mtt{i},\umtt{k}_i)$.
    Moreover, $\{\umtt{k}_1,\ldots,\umtt{k}_m\}$ forms a section relative to $[\umtt{j}]$, so that writing $s=\dimB\eta(K)$ and recalling that $\eta(K)$ is the attractor of a self-similar IFS satisfying the open set condition,
    \begin{equation}
        \sum_{i=1}^m\ctr_{\umtt{k}_i,1}^s =\ctr_{\umtt{j},1}^s.
    \end{equation}
    But $\ctr_{\mtt{i}\umtt{k}_i,1}\approx\ctr_{\mtt{i},2}$ since each $Q_i$ is an approximate square, which gives the desired result.
\end{proof}
In the next result, we provide good covers for cylinder sets using approximate squares with diameter bounded above by the height of the corresponding rectangle.
Heuristically, a cylinder set can first be decomposed into approximate squares using \cref{l:gl-pseudo-cyl-cover}, and an ``average'' approximate square itself has box dimension the same as the box dimension of $K$.
To make this notion precise, we simply reverse the order: we begin with a good cover for the box dimension of $K$, and take the image under some word $\mtt{i}$.
The image of each approximate square is a wide pseudo-cylinder, so we may apply \cref{l:gl-pseudo-cyl-cover} to complete the bound.
\begin{lemma}\label{l:gl-cyl-bound}
    Suppose $\mtt{i}\in\mathcal{I}^*$ and $0<r\leq\ctr_{\mtt{i},2}$.
    Then
    \begin{align*}
        \#\{Q\in\mathcal{S}(r):Q\subset[\mtt{i}]\}\approx \left(\frac{\ctr_{\mtt{i},2}}{r}\right)^{\dimB K}\cdot\left(\frac{\ctr_{\mtt{i},1}}{\ctr_{\mtt{i},2}}\right)^{\dimB\eta(K)}
    \end{align*}
\end{lemma}
\begin{proof}
    Fix $\mtt{i}\in\mathcal{I}^*$ and $0<r\leq\ctr_{\mtt{i},2}$.
    Write $\delta=r/\ctr_{\mtt{i},2}$, so by inspecting the proofs of \cite[Lemmas~2.1, 2.2, \& 2.3]{zbl:0757.28011}, we see that
    \begin{equation*}
        \#\mathcal{S}(\delta)\approx(1/\delta)^{\dimB K}.
    \end{equation*}
    Enumerate $\mathcal{S}(\delta)=\{Q_1,\ldots,Q_m\}$ and for each $i=1,\ldots,m$, we may write $Q_i=P(\mtt{j}_i,\umtt{k}_i)$ for some $\mtt{j}_i\in\mathcal{I}^*$ and $\umtt{k}_i\in\eta(\mathcal{I}^*)$.
    Then for each $i=1,\ldots,m$,
    \begin{equation*}
        \mathcal{Q}(\mtt{i}\mtt{j}_i,\umtt{k}_i)\subset \mathcal{S}(r)\qquad\text{and}\qquad[\mtt{i}]=\bigcup_{i=1}^m\bigcup_{Q\in\mathcal{Q}(\mtt{i}\mtt{j}_i,\umtt{k}_i)} Q.
    \end{equation*}
    Thus by \cref{l:gl-pseudo-cyl-cover} applied to each pseudo-cylinder $P(\mtt{i}\mtt{j}_i,\umtt{k}_i)$, since $Q_i$ is an approximate square and $\ctr_{\mtt{j}_i\umtt{k}_i,1}\approx \ctr_{\mtt{j}_i,2}$,
    \begin{align*}
        \#\{Q\in\mathcal{S}(r):Q\subset[\mtt{i}]\} &=\sum_{i=1}^m\#\mathcal{Q}(\mtt{i}\mtt{j}_i,\umtt{k}_i)\\
                                                   &\approx\sum_{i=1}^m\left(\frac{\ctr_{\mtt{i}\mtt{j}_i\umtt{k}_i,1}}{\ctr_{\mtt{i}\mtt{j}_i,2}}\right)^{\dimB\eta(K)}\\
                                                   &\approx\left(\frac{\ctr_{\mtt{i},2}}{r}\right)^{\dimB K}\cdot\left(\frac{\ctr_{\mtt{i},1}}{\ctr_{\mtt{i},2}}\right)^{\dimB\eta(K)}
    \end{align*}
    as claimed.
\end{proof}
To conclude our collection of preliminary lemmas, we use the Assouad dimension of the symbolic fibre $K_{\eta(\gamma)}$ to control the size of ``column sections'' of approximate squares.
We note that the word $\mtt{i}$ appears in the hypothesis but not the conclusion: this is simply to clarify the application of this lemma when it is used in \cref{p:gl-pointwise-upper}.
\begin{lemma}\label{l:gl-fibre-cover}
    Let $\epsilon>0$ and $\gamma\in\Omega$ be arbitrary.
    Suppose $k\in\N$ and $Q_k(\gamma)=P(\mtt{i},\umtt{j})$.
    Let $\mathcal{A}$ be any section of $\mathcal{I}^*$ such that $\mathcal{A}\preccurlyeq\eta^{-1}(\umtt{j})$.
    Then
    \begin{equation*}
        \sum_{\mtt{k}\in\mathcal{A}}\ctr_{\mtt{k},2}^{\dimA K_{\eta(\gamma)}+\epsilon}\lesssim_{\epsilon,\gamma} 1.
    \end{equation*}
\end{lemma}
\begin{proof}
    The assumption on the section $\mathcal{A}$ precisely means that $\{\mtt{i}\mtt{k}:\mtt{k}\in\mathcal{A}\}$ is a section relative to $\mtt{i}$ in $\mathcal{F}_{\eta(\gamma)}$.
    Then by \cref{p:Assouad-disc-packing} applied to the metric space $\Omega(\eta(\gamma))$ (recalling that $\dimA\Omega(\eta(\gamma))=\dimA K_{\eta(\gamma)}$ from \cref{c:symb-Assouad}), since $\mathcal{A}$ is a section,
    \begin{equation*}
        \sum_{\mtt{k}\in\mathcal{A}}\left(\frac{\ctr_{\mtt{i}\mtt{k},2}}{\ctr_{\mtt{i},2}}\right)^{\dimA K_{\eta(\gamma)}+\epsilon}\lesssim_{\epsilon,\gamma} 1.
    \end{equation*}
    Cancelling the $\ctr_{\mtt{i},2}$ gives the desired result.
\end{proof}
Finally, by combining the various counts that we have established earlier in this section, we are now in position to compute the upper bound for the pointwise Assouad dimension.

Let us begin with an intuitive explanation for this proof.
Since $x$ is regular, we will reduce the problem of computing covers of balls to computing covers for approximate squares.
Thus suppose we fix an approximate square $P(\mtt{i},\umtt{j})$, which is the union of cylinders $\{\mtt{i}\mtt{k}:\eta(\mtt{k})=\umtt{j}\}$.
We wish to cover this set with approximate squares in $\mathcal{S}(r)$.
There are two cases.
First, the rectangle corresponding to the cylinder $\mtt{i}\mtt{k}$ has height greater than or equal to $r$, in which case we simply keep this cylinder and obtain a good bound for the cover using \cref{l:gl-cyl-bound}: this is the family $\mathcal{A}_1$.
Otherwise, the rectangle is shorter, and we instead want to cover groups of cylinders simultaneously.
Such groups are precisely wide pseudo-cylinders corresponding to elements of $\mathcal{A}_2$ and have height $r$, which we can then cover using \cref{l:gl-pseudo-cyl-cover}.
These covers are then combined using \cref{l:gl-fibre-cover}.
\begin{proposition}\label{p:gl-pointwise-upper}
    Let $K$ be a Gatzouras--Lalley carpet and suppose $x=\pi(\gamma)\in K$.
    Then
    \begin{equation*}
        \dimA(K,x)\geq \max\{\dimB K,\dimH\eta(K)+\dimA K_{\eta(\gamma)}\}
    \end{equation*}
    with equality if $x$ is regular.
\end{proposition}
\begin{proof}
    Recalling the general lower bound proven in \cref{p:gl-pointwise-tangents}, we must show that
    \begin{equation*}
        \dimA(K,x)\leq\max\{\dimB K,\dimH\eta(K)+\dimA K_{\eta(\gamma)}\}\eqqcolon \zeta
    \end{equation*}
    when $x$ is regular.
    We obtain this bound by a direct covering argument.
    We will prove that for any $k\in\N$ and approximate square $Q_k(\gamma)=P(\mtt{i},\umtt{j})$, if $0<r\leq\ctr_{\mtt{i},2}$, then
    \begin{equation}\label{l:symb-bound}
        \#\{Q\in\mathcal{S}(r):Q\subset Q_k(\gamma)\}\lesssim\left(\frac{\ctr_{\mtt{i},2}}{r}\right)^{\zeta}.
    \end{equation}
    Assuming this, since $x$ is regular, for any ball $B(x,R)$, there is an $R'\lesssim R$ and at most two approximate squares $Q_1,Q_2\in\mathcal{S}(R')$ lying in the same column such that $B(x,R)\subset \pi(Q_1)\cup\pi(Q_2)$.
    Since $Q_1,Q_2$ lie in the same column, $Q_j=Q_{k_j}(\gamma_j)$ for some $k_j\in\N$ where $\eta(\gamma_j)=\eta(\gamma)$.
    Moreover, if $0<r\leq R$ and $Q\in\mathcal{S}(r)$ is arbitrary, then $\diam\pi(Q)\lesssim r$.
    Thus \cref{l:symb-bound} immediately gives the correct bound, up to a constant factor, for $N_r(B(x,R)\cap K)$.

    It remains to prove \cref{l:symb-bound}.
    Fix an approximate square $Q_k(\gamma)=P(\mtt{i},\umtt{j})$ and suppose $0<r\leq\ctr_{\mtt{i},2}$ is arbitrary.
    First, let
    \begin{equation*}
        \mathcal{A}_0=\eta^{-1}(\umtt{j})\wedge\mathcal{F}_{\eta(\gamma)}(r/\ctr_{\mtt{i},2})\qquad\text{and}\qquad\mathcal{A}=\{\mtt{i}\mtt{k}:\mtt{k}\in\mathcal{A}_0\}.
    \end{equation*}
    We then decompose $\mathcal{A}=\mathcal{A}_1\cup\mathcal{A}_2$, where
    \begin{equation*}
        \mathcal{A}_1=\mathcal{A}\setminus\mathcal{F}_{\eta(\gamma)}(r)\qquad\text{and}\qquad\mathcal{A}_2=\mathcal{A}\cap\mathcal{F}_{\eta(\gamma)}(r).
    \end{equation*}
    First, suppose $\mtt{i}\mtt{k}\in\mathcal{A}_1$.
    Then, by definition, $\ctr_{\mtt{i}\mtt{k},2}>r$ which, by definition of $\mathcal{A}_0$, implies that $\eta(\mtt{k})=\umtt{j}$.
    Thus by \cref{l:gl-cyl-bound} applied to the cylinder $\mtt{i}\mtt{k}$ and scale $r$, since $\dimB \eta(K)\leq\dimB K$ and $\ctr_{\mtt{i}\mtt{k},1}\approx\ctr_{\mtt{i},2}$,
    \begin{equation}\label{e:A1-bound}
        \#\{Q\in\mathcal{S}(r):Q\subset[\mtt{i}\mtt{k}]\}\approx\left(\frac{\ctr_{\mtt{i}\mtt{k},2}}{r}\right)^{\dimB K}\left(\frac{1}{\ctr_{\mtt{k},2}}\right)^{\dimB\eta(K)}.
    \end{equation}
    Otherwise, suppose $\mtt{i}\mtt{k}\in\mathcal{A}_2\subset\mathcal{F}_{\eta(\gamma)}(r)$.
    Since $\mathcal{A}_0\preccurlyeq\eta^{-1}(\umtt{j})$, there is a $\umtt{j}'$ so that $\eta(\mtt{k})\umtt{j}'=\umtt{j}$.
    Thus choice of $\umtt{j}'$ ensures that
    \begin{equation*}
        P(\mtt{i}\mtt{k},\umtt{j}')=Q_k(\gamma)\cap[\mtt{i}\mtt{k}].
    \end{equation*}
    Thus by \cref{l:gl-pseudo-cyl-cover} and since $Q_k(\gamma)=P(\mtt{i},\umtt{j})$ is an approximate square,
    \begin{equation}\label{e:A2-bound}
        \#\{Q\in\mathcal{S}(r):Q\subset Q_k(\gamma)\cap[\mtt{i}\mtt{k}]\}=\#\mathcal{Q}(\mtt{i}\mtt{k},\umtt{j}')\approx \left(\frac{1}{\ctr_{\mtt{k},2}}\right)^{\dimB\eta(K)}.
    \end{equation}
    Thus by applying \cref{e:A1-bound} and \cref{e:A2-bound} to the respective components and recalling that $\ctr_{\mtt{i}\mtt{k},2}\approx r$ whenever $\mtt{i}\mtt{k}\in\mathcal{A}_2$,
    \begin{align*}
        \#\{Q&\in\mathcal{S}(r):Q\subset Q_k(\gamma)\}\\
             &=\sum_{\mtt{i}\mtt{k}\in\mathcal{A}_1}\#\{Q\in\mathcal{S}(r):Q\subset[\mtt{i}\mtt{k}]\}+\sum_{\mtt{i}\mtt{k}\in\mathcal{A}_2}\#\{Q\in\mathcal{S}(r):Q\subset Q_k(\gamma)\cap[\mtt{i}\mtt{k}]\}\\
             &\approx\sum_{\mtt{i}\mtt{k}\in\mathcal{A}_1}\left(\frac{\ctr_{\mtt{i}\mtt{k},2}}{r}\right)^{\dimB K}\left(\frac{1}{\ctr_{\mtt{k},2}}\right)^{\dimB\eta(K)}+\sum_{\mtt{i}\mtt{k}\in\mathcal{A}_2}\left(\frac{1}{\ctr_{\mtt{k},2}}\right)^{\dimB\eta(K)}\\
             &\lesssim\sum_{\mtt{i}\mtt{k}\in\mathcal{A}_1}\left(\frac{\ctr_{\mtt{i}\mtt{k},2}}{r}\right)^{\zeta}\left(\frac{\ctr_{\mtt{i},2}}{\ctr_{\mtt{i}\mtt{k},2}}\right)^{\zeta}\ctr_{\mtt{k},2}^{\dimA K_{\eta(\gamma)}}+\sum_{\mtt{i}\mtt{k}\in\mathcal{A}_2}\left(\frac{\ctr_{\mtt{i},2}}{r}\right)^{\zeta}\ctr_{\mtt{k},2}^{\dimA K_{\eta(\gamma)}}\\
             &=\left(\frac{\ctr_{\mtt{i},2}}{r}\right)^{\zeta}\sum_{\mtt{k}\in\mathcal{A}_0}\ctr_{\mtt{k},2}^{\dimA K_{\eta(\gamma)}}\\
             &\lesssim\left(\frac{\ctr_{\mtt{i},2}}{r}\right)^{\zeta}
    \end{align*}
    where the last line follows by \cref{l:gl-fibre-cover} applied to the section $\mathcal{A}_0$.
    Thus \cref{l:symb-bound} follows, and therefore our desired result.
\end{proof}
\subsection{Dimensions of level sets of pointwise Assouad dimension}\label{ss:gl-level-set}
Given an index $i\in\mathcal{I}$, let $\Phi_{\eta(i)}$ denote the IFS corresponding to the column containing the index $i$, that is
\begin{equation*}
    \Phi_{\eta(i)}=\{S_{j,2}:j\in\mathcal{I}\text{ and }\eta(j)=\eta(i)\}.
\end{equation*}
Now given a word $\gamma=(i_n)_{n=1}^\infty\in\Omega$, recall that the symbolic slice $K_{\eta(\gamma)}$ is the non-autonomous self-similar set associated with the IFS $\{\Phi_{\eta(i_n)}\}_{n=1}^\infty$.
Since there are only finitely many choices for the $\Phi_{\eta(i_n)}$, the hypotheses of \cref{t:non-auto-Assouad} are automatically satisfied and
\begin{equation*}
    \dimA K_{\eta(\gamma)}=\lim_{m\to\infty}\sup_{n\in\N}\theta_{\eta(\gamma)}(n,m)
\end{equation*}
where
\begin{equation*}
  \sum_{(j_1,\ldots,j_m)\in\eta^{-1}(\eta(i_{n+1},\ldots,i_{n+m}))}\prod_{k=1}^m\ctr_{j_k,2}^{\theta_{\eta(\gamma)}(n,m)}=1.
\end{equation*}
We now obtain our main formula for the pointwise Assouad dimension of arbitrary points in Gatzouras--Lalley carpets.
\begin{theorem}\label{t:gl-tangent-formula}
    Let $K$ be a Gatzouras--Lalley carpet.
    Then for every $x\in K$ with $x=\pi(\gamma)$, there is an $F\in\Tan(K,x)$ with $\mathcal{H}^s(F)\gtrsim 1$ where
    \begin{align*}
        s\coloneqq{}&\dimB\eta(K)+\dimA K_{\eta(\gamma)}\\
        ={}&\dimB\eta(K)+\lim_{m\to\infty}\sup_{n\in\N}\theta_{\eta(\gamma)}(n,m)
    \end{align*}
    In particular,
    \begin{equation*}
        \max\{\dimH F:F\in\Tan(K,x)\}\geq s\quad\text{and}\quad \dimA(K,x)\geq\max\{s,\dimB K\}
    \end{equation*}
    where both inequalities are equalities if $x$ is regular.
    In particular, if $\eta(K)$ satisfies the strong separation condition then equality holds for all $x\in K$.
\end{theorem}
\begin{proof}
    By \cref{p:gl-pointwise-tangents}, there is an $F\in\Tan(K,x)$ so that
    \begin{equation*}
        \mathcal{H}^{\dimH\eta(K)+\dimA K_{\eta(\gamma)}}(F)\gtrsim 1.
    \end{equation*}
    Moreover, $\dimA K_{\eta(\gamma)}=\lim_{m\to\infty}\sup_{n\in\N}\theta_{\eta(\gamma)}(n,m)$ by \cref{t:non-auto-Assouad}.
    The formula for $\dimA(K,x)$, including the case when $x$ is regular, then follows by \cref{p:gl-pointwise-upper}.

    If $x$ is regular, it moreover follows from \cref{l:square-structure} that for any $F\in\Tan(K,x)$, there is a similarity map $h$ and a weak tangent $E\in\Tan(K_{\eta(\gamma)})$ so that $h(F)\subset\eta(K)\times E$.
    Since $\dimB\eta(K)=\dimH \eta(K)$,
    \begin{equation*}
        \dimH F=\dimH h(F)\leq\dimB \eta(K)+\dimH E\leq\dimB\eta(K)+\dimA K_{\eta(\gamma)}
    \end{equation*}
    as required.

    Finally, we recall that if $\eta(K)$ satisfies the strong separation condition, then each $x\in K$ is regular by \cref{l:gl-reg}~\cref{im:vssc}.
\end{proof}
Our next goal is to prove that the set of pointwise Assouad dimensions forms an interval.
First, for $\mtt{i}\in\mathcal{I}^n$, let $t(\mtt{i})$ be chosen so that
\begin{equation*}
    \sum_{\substack{\mtt{j}\in\mathcal{I}^n\\\eta(\mtt{j})=\eta(\mtt{i})}}\ctr_{\mtt{j},2}^{t(\mtt{i})}=1.
\end{equation*}
Equivalently, the function $t$ is chosen precisely so that
\begin{equation*}
    \theta_{\eta(\gamma)}(n,m)=t(\gamma_{n+1},\ldots,\gamma_{n+m}).
\end{equation*}
We now have the following result.
\begin{lemma}\label{l:connector-words}
    Let $K$ be a Gatzouras--Lalley carpet and suppose $\dimL K<\alpha<\dimA K$.
    Then for all $k_0\in\N$ sufficiently large, for all $n\in\N$ there is $\mtt{i}_n\in\mathcal{B}_{k_0}^n\subset\mathcal{I}^{k_0n}$ satisfying
    \begin{equation*}
        \lim_{m\to\infty}\sup_{n\in\N} t(\mtt{i}_{n+1}\cdots\mtt{i}_{n+m})=\alpha-\dimB\eta(K).
    \end{equation*}
\end{lemma}
\begin{proof}
    First, fixing any interior word $\mtt{j}\in\mathcal{I}^*$ and $i\in\mathcal{I}$ so that $\dimA K=\dimB\eta(K)+t(i)$,
    \begin{equation*}
        \dimA K=\dimB\eta(K)+\lim_{k\to\infty}t(\mtt{j}i^k);
    \end{equation*}
    and similarly for the lower dimension.
    Thus for all sufficiently large $k_0$, there are words $\mtt{j}_L,\mtt{j}_A\in\mathcal{B}_{k_0}$ so that
    \begin{equation*}
        \dimB\eta(K)+t(\mtt{j}_L)<\alpha<\dimB\eta(K)+t(\mtt{j}_A).
    \end{equation*}
    We inductively construct $(\mtt{j}_{L,k},\mtt{j}_{A,k})_{k=1}^\infty$ so that, for each $k\in\N$,
    \begin{enumerate}[nl]
        \item\label{im:jA-cnst} $\displaystyle\alpha-\dimB\eta(K)-\tfrac{1}{k}\leq t(\mtt{j}_{L,k})\leq \alpha-\dimB\eta(K)$,
        \item\label{im:jL-cnst} $\alpha-\dimB\eta(K)\leq t(\mtt{j}_{A,k})\leq \dimA K+\dimB\eta(K)+\frac{1}{k}$,
        \item\label{im:subword} $\mtt{j}_{L,k},\mtt{j}_{A,k}\in\mathcal{B}_{k_0}^*$ and, for $k\geq 2$, $\mtt{j}_{L,k},\mtt{j}_{A,k}\in\{\mtt{j}_{L,k-1},\mtt{j}_{A,k-1}\}^*$, and
        \item\label{im:long} $|\mtt{j}_{L,k}|\geq k$ and $|\mtt{j}_{A,k}|\geq k$.
    \end{enumerate}
    First, set $\mtt{j}_{L,1}=\mtt{j}_L$ and $\mtt{j}_{A,1}=\mtt{j}_A$ which clearly satisfy the desired properties.
    Now suppose we have constructed $\mtt{j}_{L,k}$ and $\mtt{j}_{A,k}$.
    Since $t(\mtt{j}_{A,k})\geq\alpha-\dimB\eta(K)$, for any $m\in\N$,
    \begin{equation*}
        \lim_{n\to\infty} t(\mtt{j}_{L,k}^m \mtt{j}_{A,k}^n)\geq\alpha-\dimB\eta(K).
    \end{equation*}
    Moreover, $t(\mtt{j}_{L,k}^m)\leq\alpha-\dimB\eta(K)$ and, by taking $m\geq k$ sufficiently large and applying \cref{t:non-auto-Assouad}~\cref{i:cont}, for all $n\in\N$ sufficiently large,
    \begin{equation*}
        |t(\mtt{j}_{L,k}^m\mtt{j}_{A,k}^{n+1})-t(\mtt{j}_{L,k}^m\mtt{j}_{A,k}^{n})|\leq\frac{1}{k+2}<\frac{1}{k+1}.
    \end{equation*}
    Combining these two observations, there is a pair $m,n$ so that $\mtt{j}_{A,k+1}\coloneqq\mtt{j}_{L,k}^m\mtt{j}_{A,k}^n\in\mathcal{B}_{k_0}^*$ satisfies conditions \cref{im:jA-cnst} and \cref{im:long}.
    The identical argument gives $\mtt{j}_{L,k+1}\in\mathcal{B}_{k_0}^*$ satisfying \cref{im:jL-cnst}, as claimed.

    To complete the proof, since $\mtt{j}_{L,k}\in\mathcal{B}_{k_0}^*$ for all $k\in\N$, we may identify the sequence $(\mtt{j}_{L,k})_{k=1}^\infty$ with a sequence $(\mtt{i}_n)_{n=1}^\infty$ where $\mtt{i}_n\in\mathcal{B}_{k_0}$ for all $n\in\N$.
    It immediately follows from \cref{im:jA-cnst} and \cref{im:long} that
    \begin{equation*}
        \lim_{m\to\infty}\sup_{n\in\N} t(\mtt{i}_{n+1}\cdots\mtt{i}_{n+m})\geq \alpha-\dimB\eta(K).
    \end{equation*}
    To establish the converse bound, it suffices to show for every $k\in\N$ that
    \begin{equation*}
        \lim_{m\to\infty}\sup_{n\in\N} t(\mtt{i}_{n+1}\cdots\mtt{i}_{n+m})\leq\alpha-\dimB\eta(K)+\frac{1}{k}.
    \end{equation*}
    By \cref{im:subword}, for all $k\in\N$, there is a $K\in\N$ so that for all $n\geq K$, $\mtt{i}_n\in\{\mtt{j}_{L,k},\mtt{j}_{A,k}\}^*$.
    For each $\ell\in\N$, write $\mtt{k}_\ell=\mtt{i}_{K\ell+1}\cdots\mtt{i}_{K(\ell+1)}$ and note that $\mtt{k}_\ell\in\{\mtt{j}_{L,k},\mtt{j}_{A,k}\}^*$ for all $\ell\in\N$.
    Thus for any $n,m\in\N$,
    \begin{equation*}
        t(\mtt{k}_{\ell+1}\cdots \mtt{k}_{\ell+m})\leq \frac{1}{m}\sum_{i=1}^m t(\mtt{k}_{\ell+i})\leq \alpha-\dimB\eta(K)+\frac{1}{k}.
    \end{equation*}
    But by the property of $t$ established in \cref{t:non-auto-Assouad}~\cref{e:sub-lim},
    \begin{equation*}
        \lim_{m\to\infty}\sup_{n\in\N} t(\mtt{i}_{n+1}\cdots\mtt{i}_{n+m})
        =\lim_{m\to\infty}\sup_{n\in\N} t(\mtt{k}_{n+1}\cdots\mtt{k}_{n+m})
    \end{equation*}
    which gives the claim.
\end{proof}
To conclude this section, we assemble the results proven in the prior two sections to obtain our main result.
\begin{theorem}\label{t:gl-large-Assouad}
    Let $K$ be a Gatzouras--Lalley carpet.
    Then for any $\dimB K\leq\alpha\leq\dimA K$,
    \begin{equation}\label{e:full-off}
        \dimH\{x\in K:\dimA(K,x)=\alpha\}=\dimH K.
    \end{equation}
    Otherwise, if $\alpha\notin[\dimB K,\dimA K]$, then $\{x\in K:\dimA(K,x)=\alpha\}=\varnothing$.
    However,
    \begin{equation}\label{e:meas-zero}
        \mathcal{H}^{\dimH K}\bigl(\{x\in K:\dimA(K,x)\neq\dimA K\}\bigr)=0.
    \end{equation}
\end{theorem}
\begin{proof}
    Note that if $\dimB K=\dimA K$, then $\dimA(K,x)=\dimA K$ for all $x\in K$ and the results are clearly true.
    Thus we may assume that $\dimH K< \dimB K<\dimA K$.

    We first establish \cref{e:full-off}.
    Let $\epsilon>0$ be arbitrary and $\dimB K\leq\alpha\leq\dimA K$.
    Apply \cref{p:gl-large-subsystem} and get $k\in\N$ and a family $\mathcal{J}\subset\mathcal{B}_k$ with corresponding attractor $K_\epsilon$ satisfying $\dimH K-\epsilon\leq\dimH K_\epsilon=\dimA K_\epsilon$ and $\dimB\eta(K)-\epsilon\leq\dimB\eta(K_\epsilon)$.
    If $\alpha<\dimA K$, iterating the system if necessary, by \cref{l:connector-words} get a sequence $(\mtt{i}_n)_{n=1}^\infty$ with $\mtt{i}_n\in\mathcal{B}_k$ for all $n\in\N$ and moreover
    \begin{equation}\label{e:large-subword}
        \lim_{m\to\infty}\sup_{n\in\N}t(\mtt{i}_{n+1}\cdots\mtt{i}_{n+m})=\alpha-\dimB\eta(K).
    \end{equation}
    If instead $\alpha=\dimA K$, instead simply take $\mtt{i}_n= i_0^k$ where $i_0\in\mathcal{I}$ is any word such that $\dimA K=\dimB\eta(K)+t(i_0)$.
    Note that $t(\mtt{j})=\dimA K_\epsilon-\dimB\eta(K_\epsilon)$ for any $\mtt{j}\in\mathcal{J}$.
    Thus by taking $\epsilon$ to be sufficiently small, we may assume that $t(\mtt{j})\leq\alpha-\dimB\eta(K)$ for all $\mtt{j}\in\mathcal{J}$.

    Now, let $(N_k)_{k=1}^\infty$ be a sequence of natural numbers satisfying $\lim_{k\to\infty}N_k/k=\infty$ and write
    \begin{equation*}
        \Omega_0=\prod_{k=1}^\infty\mathcal{J}^{N_k}\times\{\mtt{i}_1\}\times\cdots\times\{\mtt{i}_k\}.
    \end{equation*}
    By taking each $N_k$ to be sufficiently large, we may ensure that $\dimH \pi(\Omega_0)\geq\dimH K_\epsilon-\epsilon$.
    Fix $\gamma\in\Omega_0$: it remains to verify that $\dimA(K,\pi(\gamma))=\alpha$.
    Since $\gamma\in\mathcal{B}_k^{\N}$, $\pi(\gamma)$ is a regular point of $K$ by \cref{l:gl-reg}~\cref{im:all-subword}.
    By passing to the subsystem induced by $\mathcal{B}_k\subset\mathcal{I}^k$, write $\gamma=(\mtt{k}_k)_{k=1}^\infty$ where $\mtt{k}_k\in\mathcal{B}_k$.
    Thus by \cref{t:gl-tangent-formula,t:non-auto-Assouad},
    \begin{equation*}
        \dimA(K,x)=\max\bigl\{\dimB K, \lim_{m\to\infty}\sup_{n\in\N}t(\mtt{k}_{n+1}\cdots \mtt{k}_{n+m})+\dimB\eta(K)\bigr\}.
    \end{equation*}
    Since $\mtt{i}_1\cdots\mtt{i}_m$ appears as a subword of $\gamma$ for arbitrarily large $m$, by \cref{e:large-subword} and since $\alpha>\dimB K$, it follows that $\dimA(K,x)\geq\alpha$.

    We now obtain the upper bound.
    Let $\epsilon>0$ be arbitrary.
    By \cref{e:large-subword}, there is an $\ell_0\in\N$ so that whenever $\ell\geq\ell_0$, we have $t(\mtt{i}_{j+1}\cdots\mtt{i}_{j+\ell})\leq \alpha-\dimB\eta(K)+\epsilon$.
    Let $C$ be the implicit constant from \cref{t:non-auto-Assouad}~\cref{i:cont} and let $m$ be sufficiently large so that $C\ell_0/m\leq\epsilon$.
    Since $\lim_{k\to\infty}N_k/k=\infty$, for all $n$ sufficiently large, there is a $j\in\N$ so that
    \begin{equation*}
        \mtt{k}_{n+1}\cdots \mtt{k}_{n+m}=\mtt{j}_1\cdots\mtt{j}_{m-\ell}\mtt{i}_{j+1}\cdots\mtt{i}_{j+\ell}.
    \end{equation*}
    Thus for $m,n$ sufficiently large, if $\ell\geq \ell_0$, by \cref{t:non-auto-Assouad},
    \begin{align*}
        t(\mtt{k}_{n+1}\cdots \mtt{k}_{n+m})&\leq \max\bigl\{t(\mtt{j}_1\cdots\mtt{j}_{m-\ell}),t(\mtt{i}_{j+1}\cdots\mtt{i}_{j+\ell})\}\\
                                            &\leq \alpha-\dimB\eta(K)+\epsilon
    \end{align*}
    and similarly if $\ell<\ell_0$, recalling that $t(\mtt{i}_{j+1}\cdots\mtt{i}_{j+\ell})\leq 1$, by \cref{t:non-auto-Assouad} recalling the definition of $C$,
    \begin{equation*}
        t(\mtt{k}_{n+1}\cdots \mtt{k}_{n+m})\leq\alpha-\dimB\eta(K)+\epsilon.
    \end{equation*}
    Therefore
    \begin{equation*}
        \limsup_{m\to\infty}\limsup_{n\to\infty}t(\mtt{k}_{n+1}\cdots \mtt{k}_{n+m})\leq \alpha-\dimB\eta(K)+\epsilon
    \end{equation*}
    and since $\epsilon>0$ was arbitrary,
    \begin{equation*}
        \lim_{m\to\infty}\sup_{n\in\N}t(\mtt{k}_{n+1}\cdots \mtt{k}_{n+m})=\limsup_{m\to\infty}\limsup_{n\to\infty}t(\mtt{k}_{n+1}\cdots \mtt{k}_{n+m})\leq \alpha-\dimB\eta(K)
    \end{equation*}
    so that $\dimA(K,x)\leq\alpha$, as claimed.
    Of course, we recall as well that $\dimB K\leq\dimA(K,x)\leq\dimA K$ by \cref{p:tan-lower}.

    We finally consider the points $x$ such that $\dimA(K,x)<\dimA K$.
    Let $i_0\in\mathcal{I}$ be such that $\dimA K=\dimB\eta(K)+t(i_0)$.
    Let
    \begin{equation*}
        \mathcal{J}_M\coloneqq \{(i_1,\ldots,i_M)\in\mathcal{I}^M:(i_1,\ldots,i_M)\neq(i_0,\ldots,i_0)\}
    \end{equation*}
    have attractor $K_M\subset K$.
    Since $\mathcal{J}_M$ is a proper subsystem, $\dimH K_M<\dimH K$ so that $\mathcal{H}^{\dimH K}(K_M)=0$.
    Now let $x\in K$ have $\dimA(K,x)<\dimA K$.
    Suppose $x=\pi(\gamma)$ where $\gamma=(i_n)_{n=1}^\infty$, so that
    \begin{equation*}
        \dimA(K,x)\geq \max\left\{\dimB K,\dimB\eta(K)+\lim_{m\to\infty}\sup_{n\in\N}t(i_{n+1},\ldots,i_{n+m})\right\}.
    \end{equation*}
    Since $\dimA(K,x)<\dimA K$,
    \begin{equation*}
        \lim_{m\to\infty}\sup_{n\in\N}t(i_{n+1},\ldots,i_{n+m})<t(i_0).
    \end{equation*}
    In particular, there is a constant $M$ so that $\gamma$ does not contain $i_0^M$ as a subword.
    Thus $x\in K_M$ for some $M$ and therefore
    \begin{equation*}
        \mathcal{H}^{\dimH K}\left(\{x\in K:\dimA(K,x)<\dimA K\}\right)\leq\sum_{M=1}^\infty\mathcal{H}^{\dimH K}(K_M)=0
    \end{equation*}
    as required.
\end{proof}
\begin{remark}
    We recall that if $K$ is a Gatzouras--Lalley carpet, then $\mathcal{H}^{\dimH K}(K)>0$, with $\mathcal{H}^{\dimH K}(K)<\infty$ if and only if $K$ is Ahlfors regular; see \cite{zbl:0757.28011}.
    In particular, the positivity of the Hausdorff measure guarantees that the claim \cref{e:meas-zero} in \cref{t:gl-large-Assouad} is not vacuous; and, if the Hausdorff measure is finite, \cref{t:gl-large-Assouad} is trivial.
\end{remark}

\section{Tangent structure and dimension of Barański carpets} \label{sec:baranski}
In this section, we adapt the strategy of \cref{s:carpet-tangents} from Gatzouras--Lalley to Barański carpets, allowing the dominant direction of contraction to change along the orbit.
As before, we express the Hausdorff and Assouad dimensions in terms of suitable symbolic slices and metric trees of approximate squares, but now in both horizontal and vertical directions and with an additional decomposition of the coding space according to the (local) direction of maximal contraction.
The lower bounds using the non-autonomous fibre theory from \cref{sec:non-autonomous} are morally quite similar to the Gatzouras--Lalley case.
In contrast, the upper bounds require more work: even though the fibre may be globally contracting, this may not be the case uniformly over all scales.
This is handled in \cref{ss:bar-up} and constitutes the majority of the work in this section.

\subsection{Dimensions and decompositions of Barański carpets} \label{sec:baranski-dim}
Recall the definition of the Barański carpet and basic notation from \cref{ss:carpet-def}.
Suppose $K$ is a Barański carpet and $\gamma\in\Omega$ is arbitrary.
For each $k\in\N$, we define a probability vector $\bm{\xi}_k(\gamma)$ by the rule
\begin{equation*}
    \bm{\xi}_k(\gamma)_i=\frac{\#\{1\leq \ell\leq k:\gamma_\ell=i\}}{k}\quad\text{for each } i\in\mathcal{I}.
\end{equation*}
In other words, $\bm{\xi}_k(\gamma)$ is the distribution of the letter frequencies in the first $k$ letters of $\gamma$.
We then define
\begin{equation*}
    \Gamma_k(\gamma)=\frac{\chi_1(\bm{\xi}_k(\gamma))}{\chi_2(\bm{\xi}_k(\gamma))}.
\end{equation*}
The function $\Gamma_k$ induces a partition $\Omega=\Omega_0\cup\Omega_1\cup\Omega_2$ by
\begin{align*}
    \Omega_0&=\{\gamma:\liminf_{k\to\infty}\Gamma_k(\gamma)\leq 1\leq \limsup_{k\to\infty}\Gamma_k(\gamma)\}\\
    \Omega_1&=\{\gamma:\limsup_{k\to\infty}\Gamma_k(\gamma)<1\}\\
    \Omega_2&=\{\gamma:1<\liminf_{k\to\infty}\Gamma_k(\gamma)\}.
\end{align*}
We now recall the dimensional formula for a general Barański carpet.
First, we decompose $\mathcal{P}=\mathcal{P}_1\cup\mathcal{P}_2$ where
\begin{equation*}
    \mathcal{P}_j=\{\bm{w}\in\mathcal{P}:\chi_j(\bm{w})\leq\chi_{j'}(\bm{w})\}.
\end{equation*}
Now given a measure $\bm{w}\in\mathcal{P}_j$, recall \cite[Corollary~5.2]{zbl:1116.28008} which states that
\begin{equation*}
    \dimH \pi_*\bm{w}^{\N}=\frac{H(\eta_j(\bm{w}))}{\chi_j(\bm{w})}+\frac{H(\bm{w})-H(\eta_j(\bm{w}))}{\chi_{j'}(\bm{w})}.
\end{equation*}
Here and for the remainder of this document, for notational simplicity, given $j=1$ we write $j'=2$ and given $j=2$ we write $j'=1$.

We also introduce some notation for symbolic slices both in the horizontal and vertical directions.
Given $\gamma=(i_n)_{n=1}^\infty\in\Omega$ and $j\in\{1,2\}$, let $\theta_{\eta_j(\gamma),j}$ be defined by the rule
\begin{equation*}
  \sum_{(j_1,\ldots,j_m)\in\eta_j^{-1}(\eta_j(i_{n+1},\ldots,i_{n+m}))}\prod_{k=1}^m\ctr_{j_k,j'}^{\theta_{\eta_j(\gamma),j}(n,m)}=1.
\end{equation*}
The value $\theta_{\eta(\gamma)}=\theta_{\eta_1(\gamma),1}$ was defined previously in the context of a Gatzouras--Lalley carpet.
As is the case with a Gatzouras--Lalley carpet, if we denote by $K_{\eta_j(\gamma),j}$ the non-autonomous self-similar set associated with the non-autonomous self-similar IFS $\{S_{i,j'}:i\in\eta_j^{-1}(\eta_j(\gamma_k))\}_{k=1}^\infty$, then
\begin{equation*}
    \dimA K_{\eta_j(\gamma),j}=\lim_{m\to\infty}\sup_{n\in\N}\theta_{\eta_j(\gamma),j}(n,m).
\end{equation*}
Assuming $\eta_1(K)$ (resp.~$\eta_2(K)$) satisfies the SSC, then $K_{\eta_1(\gamma),1}$ (resp.~$K_{\eta_2(\gamma),2}$) is precisely the intersection of $K$ with the vertical (resp.~horizontal) line containing $x=\pi(\gamma)$.

We now recall \cite[Theorem~2.12]{zbl:1305.28021} concerning the Assouad dimension and \cite[Theorems~A and~B]{zbl:1116.28008} on the Hausdorff and box dimensions of Barański carpets.
While these results are not stated explicitly in this form, the relevant details can be obtained directly by inspecting the proofs.
\begin{proposition}[\cite{zbl:1305.28021,zbl:1116.28008}]\label{p:bar-dims}
    Let $K$ be a Barański carpet such that $\Omega_1\neq\varnothing$ and $\Omega_2\neq\varnothing$.
    Then:
    \begin{enumerate}[nl,r]
        \item\label{im:bar-adim} For each $j=1,2$,
            \begin{equation*}
                \dimH\pi(\Omega_0\cup\Omega_j)\leq d_j
            \end{equation*}
            where
            \begin{equation*}
                d_j=\max_{\bm{w}\in\mathcal{P}_j}\left(\frac{H(\eta_j(\bm{w}))}{\chi_j(\bm{w})}+\frac{H(\bm{w})-H(\eta_j(\bm{w}))}{\chi_{j'}(\bm{w})}\right).
            \end{equation*}
            In particular, $\dimH K=\max\{d_1,d_2\}$.
        \item\label{it:bar-bdim} We have
        \begin{equation*}
          \dimB K=\max_{j=1,2}\{\dimB\eta_j(K)+z_j\}
        \end{equation*}
        where $z_j$ is the unique solution to the equation
        \begin{equation*}
          \sum_{i\in\mathcal{I}}\ctr_{i,j}^{\dimB\eta_j(K)}\ctr_{i,j'}^{z_j}=1.
        \end{equation*}
        \item\label{im:bar-hdim} We have
            \begin{equation*}
                \dimA K=\max_{j=1,2}\left\{\dimB\eta_j(K)+t_j\right\}
            \end{equation*}
            where
            \begin{equation*}
                t_j=\max_{\underline{\ell}\in\eta_j(\mathcal{I})}t_j(\underline{\ell})
            \end{equation*}
            and $t_j(\underline{\ell})$ is the unique solution to the equation
            \begin{equation*}
              \sum_{i\in\eta_j^{-1}(\underline{\ell})}\ctr_{i,j'}^{t_j(\underline{\ell})}=1.
            \end{equation*}
    \end{enumerate}
\end{proposition}

\subsection{Pointwise Assouad dimension along uniformly contracting sequences}\label{ss:bar-up}
In this section, we state a generalization of our results on Gatzouras--Lalley carpets to Barański carpets, with the caveat that we restrict our attention to points coded by sequences which contract uniformly in one direction.
The approach is similar to the Gatzouras--Lalley case so we only include the details when the proofs diverge.
Handling more general sequences would result in a more complicated formula for the pointwise Assouad dimension depending on the scales at which the contraction ratio is greater in one direction than the other, which we will not treat here.

We begin by defining the analogues of pseudo-cylinders and approximate squares.
Fix $j=1,2$.
Suppose $\mtt{i}\in\mathcal{I}^k$ and $\umtt{j}\in\eta_j(\mathcal{I}^\ell)$.
We then write
\begin{equation*}
    P_j(\mtt{i},\umtt{j})=\{\gamma=(i_n)_{n=1}^\infty\in\Omega:(i_1,\ldots,i_k)=\mtt{i}\text{ and }\eta_j(i_{k+1},\ldots,i_{k+l})=\umtt{j}\}.
\end{equation*}
Now let $\gamma\in\Omega$ be arbitrary and let $k\in\N$.
Let $j$ be chosen so that $\ctr_{\gamma\npre{k},{j}}\geq \ctr_{\gamma\npre{k},j'}$.
We then let $L_{k}(\gamma)\geq k$ be the minimal integer so that
\begin{equation*}
    \ctr_{\gamma\npre{L_{k,j}(\gamma)},{j}}<\ctr_{\gamma\npre{k},j'}.
\end{equation*}
Write $\gamma\npre{L_{k,j}(\gamma)}=\mtt{i}\mtt{j}$ and define the approximate square
\begin{equation*}
    Q_{k}(\gamma) = P_j(\mtt{i},\eta_j(\mtt{j})).
\end{equation*}
Finally, we call a pseudo-cylinder \emph{wide} if $P_j(\mtt{i},\umtt{j})$ contains an approximate square $P_j(\mtt{i},\umtt{k})$; otherwise, we call the pseudo-cylinder \emph{tall}.

In the case when the Barański carpet is in fact a Gatzouras--Lalley carpet, these definitions with $j=1$ coincide with the definitions in the Gatzouras--Lalley case.

Next, the collection of approximate squares forms a metric tree when equipped with the valuation $\rho(P_j(\mtt{i},\eta_j(\mtt{j})))=\ctr_{\mtt{i},j'}$.
Note that for each approximate square $Q$, there is a unique choice for $j$ except precisely when $\ctr_{\gamma\npre{k},j}=\ctr_{\gamma\npre{k},j'}$, so indeed $\rho$ is well-defined.

Similarly as in the Gatzouras--Lalley case, given a pseudo-cylinder $P_j(\mtt{i},\umtt{j})$, we write
\begin{equation*}
    \mathcal{Q}_j(\mtt{i},\umtt{j})=\max\{\mathcal{A}:\mathcal{A}\text{ is a section of }\mathcal{S}\text{ relative to }P_j(\mtt{i},\umtt{j})\}
\end{equation*}
where $\mathcal{S}$ is the collection of all approximate squares and the maximum is with respect to the partial ordering on sections.
That the maximum always exists follows from the properties of the meet.
In the case when the pseudo-cylinder is wide, this coincides precisely with the definition in the Gatzouras--Lalley case.

However, unlike in the Gatzouras--Lalley case, we will also have to handle tall pseudo-cylinders, which have a more complex structure.
This additional structure is handled in the following covering lemma.
\begin{lemma}\label{l:bar-covers}
    \begin{enumerate}[nl,r]
        \item\label{im:bar-wide-pseudo} Let $P_j(\mtt{i},\umtt{j})$ be a wide pseudo-cylinder.
            Then
            \begin{equation*}
                \#\mathcal{Q}_j(\mtt{i},\umtt{j})\approx\left(\frac{\ctr_{\mtt{i}\umtt{j},j}}{\ctr_{\mtt{i},j'}}\right)^{\dimB\eta_j(K)}.
            \end{equation*}
        \item\label{im:bar-tall-pseudo} Let $P_j(\mtt{i},\umtt{j})$ be a tall pseudo-cylinder.
            Then
            \begin{equation*}
                \#\mathcal{Q}_j(\mtt{i},\umtt{j})\lesssim\left(\frac{\ctr_{\mtt{i},j'}}{\ctr_{\mtt{i}\umtt{j},j}}\right)^{\dimB\eta_{j'}(K)}.
            \end{equation*}
        \item\label{im:bar-cyl-cover} Let $\epsilon>0$ be arbitrary.
            Suppose $\mtt{i}\in\mathcal{I}^*$ and let $j$ be chosen so that $\ctr_{\mtt{i},j'}\leq\ctr_{\mtt{i},j}$.
            Let $0<r\leq\ctr_{\mtt{i},j}$.
            Then
            \begin{equation*}
                \#\{Q\in\mathcal{S}(r):Q\subset[\mtt{i}]\}\lesssim_\epsilon\left(\frac{\ctr_{\mtt{i},j'}}{r}\right)^{\dimB K+\epsilon}\cdot\left(\frac{\ctr_{\mtt{i},j}}{\ctr_{\mtt{i},j'}}\right)^{\dimB\eta_j(K)}.
            \end{equation*}
        \item\label{im:bar-fibre-cover} Let $\epsilon>0$ and $\gamma\in\Omega$ be arbitrary.
            Suppose $k\in\N$ and $j=1,2$ are such that $Q_k(\gamma)=P_j(\mtt{i},\umtt{j})$.
            Let $\mathcal{A}$ be any section of $\mathcal{I}^*$ satisfying $\mathcal{A}\preccurlyeq\eta_j^{-1}(\umtt{j})$.
            Then
            \begin{equation*}
                \sum_{\mtt{k}\in\mathcal{A}}\ctr_{\mtt{k},j'}^{\dimA K_{\eta_j(\gamma),j}+\epsilon}\lesssim_{\epsilon,\gamma} 1.
            \end{equation*}
    \end{enumerate}
\end{lemma}
\begin{proof}
    The proof of \cref{im:bar-wide-pseudo} is identical to the proof given in \cref{l:gl-pseudo-cyl-cover} and similarly the proof of \cref{im:bar-fibre-cover} is identical to that of \cref{l:gl-fibre-cover}.

    We now prove \cref{im:bar-tall-pseudo}.
    In order to do this, we must understand the structure of the pseudo-cylinder $P_j(\mtt{i},\umtt{j})$.
    Heuristically, when (for instance) $j=1$, $P_j(\mtt{i},\umtt{j})$ is a union of cylinders which fall into one of two types: those which are tall, and those which are wide.
    If a cylinder is tall, we apply \cref{im:bar-wide-pseudo} in the opposite direction to cover it with approximate squares, and if a cylinder is wide, we group nearby cylinders together to form approximate squares.
    We then combine these counts using the slice dimension $t_j$, which is bounded above by $\dimB\eta_{j'}(K)$.

    Write $\mathcal{A}=\eta_j^{-1}(\umtt{j})$ and partition $\mathcal{A}=\mathcal{A}_1\cup\mathcal{A}_2$ where
    \begin{equation*}
        \mathcal{A}_1=\{\mtt{k}\in\mathcal{A}:\ctr_{\mtt{i}\mtt{k},j'}\geq\ctr_{\mtt{i}\umtt{j},j}\}\qquad\text{and}\qquad\mathcal{A}_2=\mathcal{A}\setminus\mathcal{A}_1.
    \end{equation*}
    First, for $\mtt{k}\in\mathcal{A}_1$, note that $P_{j'}(\mtt{i}\mtt{k},\varnothing)$ is a wide pseudo-cylinder and we set
    \begin{equation*}
        \mathcal{B}_1=\bigcup_{\mtt{k}\in\mathcal{A}_1}\mathcal{Q}_{j'}(\mtt{i}\mtt{k},\varnothing).
    \end{equation*}
    By applying \cref{im:bar-wide-pseudo}, since $\ctr_{\mtt{i}\mtt{k},j}\approx \ctr_{\mtt{i}\umtt{j},j}$,
    \begin{equation}\label{e:B1-bound}
        \#\mathcal{B}_1=\sum_{\mtt{k}\in\mathcal{A}_1}\#\mathcal{Q}_{j'}(\mtt{i}\mtt{k},\varnothing)\approx\sum_{\mtt{k}\in\mathcal{A}_1}\left(\frac{\ctr_{\mtt{i}\mtt{k},j'}}{\ctr_{\mtt{i}\umtt{j},j}}\right)^{\dimB\eta_{j'}(K)}
    \end{equation}
    Otherwise if $\mtt{k}\in\mathcal{A}_2$, let $\mtt{l}_1(\mtt{k})$ denote the prefix of $\mtt{k}$ of maximal length so that $\ctr_{\mtt{i}\mtt{l}_1(\mtt{k}),j'}\geq\ctr_{\mtt{i}\umtt{j},j}$.
    Writing $\mtt{k}=\mtt{l}_1(\mtt{k})\mtt{l}_2(\mtt{k})$, this choice guarantees that
    \begin{equation*}
        \mathcal{B}(\mtt{k})\coloneqq P_j(\mtt{i}\mtt{l}_1(\mtt{k}),\eta_j(\mtt{l}_2(\mtt{k})))
    \end{equation*}
    is the unique approximate square contained in $[\mtt{i}]$ containing $[\mtt{i}\mtt{k}]$.
    Finally, let
    \begin{equation*}
        \mathcal{A}_2'=\{\mtt{l}_1(\mtt{k}):\mtt{k}\in\mathcal{A}_2\}\qquad\text{and}\qquad\mathcal{B}_2=\{\mathcal{B}(\mtt{k}):\mtt{k}\in\mathcal{A}_2\}.
    \end{equation*}
    We then note that, since $\ctr_{\mtt{i}\mtt{l},j'}\approx\ctr_{\mtt{i}\umtt{j},j}$ by the choice of $\mtt{l}_1(\mtt{k})$,
    \begin{equation}\label{e:B2-bound}
        \#\mathcal{B}_2\approx\sum_{\mtt{l}\in\mathcal{A}_2'}\left(\frac{\ctr_{\mtt{i}\mtt{l},j'}}{\ctr_{\mtt{i}\umtt{j},j}}\right)^{\dimB\eta_{j'}(K)}
    \end{equation}
    To conclude, observe that $\mathcal{Q}_j(\mtt{i},\umtt{j})=\mathcal{B}_1\cup\mathcal{B}_2$ and applying \cref{e:B1-bound} and \cref{e:B2-bound},
    \begin{align*}
        \#\mathcal{Q}_j(\mtt{i},\umtt{j}) &= \#\mathcal{B}_1+\#\mathcal{B}_2\\
                                          &\lesssim\sum_{\mtt{k}\in\mathcal{A}_1}\left(\frac{\ctr_{\mtt{i}\mtt{k},j'}}{\ctr_{\mtt{i}\umtt{j},j}}\right)^{\dimB\eta_{j'}(K)}+\sum_{\mtt{l}\in\mathcal{A}_2'}\left(\frac{\ctr_{\mtt{i}\mtt{l},j'}}{\ctr_{\mtt{i}\umtt{j},j}}\right)^{\dimB\eta_{j'}(K)}\\
                                          &=\left(\frac{\ctr_{\mtt{i},j'}}{\ctr_{\mtt{i}\umtt{j},j}}\right)^{\dimB\eta_{j'}(K)}\sum_{\mtt{k}\in\mathcal{A}_1\cup\mathcal{A}_2'}\ctr_{\mtt{k},j'}^{\dimB\eta_{j'}(K)}\\
                                          &\leq\left(\frac{\ctr_{\mtt{i},j'}}{\ctr_{\mtt{i}\umtt{j},j}}\right)^{\dimB\eta_{j'}(K)}
    \end{align*}
    where the last line follows since $\mathcal{A}_1\cup\mathcal{A}_2'\preccurlyeq\eta_j^{-1}(\umtt{j})$ is a section and $\dimB\eta_{j'}(K)\geq t_j(\umtt{j})$ where
    \begin{equation*}
        \sum_{\mtt{k}\in\mathcal{A}_1\cup\mathcal{A}_2'}\ctr_{\mtt{k},j'}^{t_j(\umtt{j})}=1.
    \end{equation*}
    Finally, we combine the bounds given in \cref{im:bar-wide-pseudo} and \cref{im:bar-tall-pseudo} with a similar argument to the proof of \cref{l:gl-cyl-bound} to obtain \cref{im:bar-cyl-cover}.
    Let $\epsilon>0$ be arbitrary and fix $\mtt{i}\in\mathcal{I}^*$ and $j=0,1$ so that $0<r\leq\ctr_{\mtt{i},j'}\leq\ctr_{\mtt{i},j}$.
    Write $\delta=r/\ctr_{\mtt{i},j'}$ so, recalling the proof of \cite[Theorem~B]{zbl:1116.28008},
    \begin{equation*}
        \#\mathcal{S}(\delta)\lesssim_\epsilon(1/\delta)^{\dimB K+\epsilon}.
    \end{equation*}
    Now enumerate
    \begin{equation*}
        \mathcal{S}(\delta)=\{Q_{1,j},\ldots,Q_{m_j,j}\}\cup\{Q_{1,j'},\ldots,Q_{m_{j'},j'}\}
    \end{equation*}
    where for each $z=j,j'$ and $1\leq i\leq m_z$,
    \begin{equation*}
        Q_{i,z}=P_z(\mtt{j}_{i,z},\umtt{k}_{i,z})
    \end{equation*}
    for some $\mtt{j}_{i,z}\in\mathcal{I}^*$ and $\umtt{k}_{i,z}\in\eta_z(\mathcal{I}^*)$.
    Observe that each $P_z(\mtt{i}\mtt{j}_{i,z},\umtt{k}_{i,z})$ is a wide pseudo-cylinder if $z=j$ and a tall pseudo-cylinder if $z=j'$.
    Thus we may complete the proof in the same way as \cref{l:gl-cyl-bound}, by applying \cref{im:bar-wide-pseudo} to the wide pseudo-cylinders and \cref{im:bar-tall-pseudo} to the tall pseudo-cylinders.
\end{proof}
We can now prove the following formulas for the pointwise Assouad dimension.
\begin{proposition}\label{p:baranski-tans}
    Let $K$ be a Barański carpet.
    Then for each $j=1,2$, if $\eta_j(K)$ satisfies the SSC, for all $\gamma\in\Omega_j$ and $x=\pi(\gamma)$,
    \begin{equation*}
        \dimA(K,x)=\max\{\dimB K, \dimB\eta_j(K)+\dimA K_{\eta_j(\gamma),j}\}
    \end{equation*}
    and
    \begin{equation*}
        \max\{\dimH F:F\in\Tan(K,x)\}=\dimB\eta_j(K)+\dimA K_{\eta_j(\gamma),j}.
    \end{equation*}
    Furthermore,
    \begin{equation*}
        \dimA K_{\eta_j(\gamma),j}=\lim_{m\to\infty}\sup_{n\in\N}\theta_{\eta_j(\gamma),j}(n,m)\leq \max_{\underline{\ell}\in\eta_j(\mathcal{I})} t_j(\underline{\ell}).
    \end{equation*}
\end{proposition}
\begin{proof}
    If $\gamma\in\Omega_j$, by definition there is a constant $\kappa\in(0,1)$ so that
    \begin{equation*}
        \frac{\ctr_{\gamma\npre{k},{j'}}}{\ctr_{\gamma\npre{k},j}}\lesssim\kappa^n.
    \end{equation*}
    In particular, there is a constant $\kappa'\in(0,1)$ so that each maximal cylinder $[\mtt{i}]$ contained in $Q_k(\gamma)$ has $\ctr_{\mtt{i},j'}/\ctr_{\mtt{i},j}\lesssim(\kappa')^k$, which converges to zero.
    Thus the same proof as given in \cref{p:gl-pointwise-upper} but instead applying \cref{l:bar-covers} in place of the analogous bounds for Gatzouras--Lalley carpets gives that
    \begin{equation*}
        \dimA(K,x)\leq\max\{\dimB K, \dimB\eta_j(K)+\dimA K_{\eta_j(\gamma),j}\}.
    \end{equation*}
    Similarly, the same proof as \cref{p:gl-pointwise-tangents} shows that
    \begin{equation*}
        \max\{\dimH F:F\in\Tan(K,x)\}=\dimB\eta_j(K)+\dimA K_{\eta_j(\gamma),j}.
    \end{equation*}
    Finally, using \cref{t:non-auto-Assouad},
    \begin{equation*}
        \lim_{m\to\infty}\sup_{n\in\N}\theta_{\eta_j(\gamma),j}(n,m)\leq \max_{\underline{\ell}\in\eta_j(\mathcal{I})} t_j(\underline{\ell}).
    \end{equation*}
    as required.
\end{proof}

\subsection{Barański carpets with few large tangents}
In contrast to Gatzouras--Lalley carpets, the analogue of \cref{t:gl-large-Assouad} need not hold for Barański carpets.
We first give a precise characterization of when a Barański carpet has few large tangents.
Fix the definitions of $t_j$ and $d_j$ from \cref{p:bar-dims}.

We begin with an auxiliary lemma which gives conditions for which the box dimension of a Barański carpet is strictly smaller than its Assouad dimension.
\begin{lemma} \label{lem:bar-box}
    Let $K$ be a Barański carpet such that $\Omega_1\neq\varnothing$ and $\Omega_2\neq\varnothing$.
    Suppose for one of $j=1,2$, $d_j<d_{j'}$ and $\dimB\eta_j(K)+t_j>\dimB\eta_{j'}(K)+t_{j'}$.
    Then $\dimB K<\dimA K$.
\end{lemma}
\begin{proof}
    For $j=1,2$ write $s_j=\dimB\eta_j(K)$ and recall the number $z_j$ from \cref{p:bar-dims}~\cref{it:bar-bdim}.
    Recall that $s_j$ is the unique solution to the equation $\sum_{\underline{\ell}\in\eta_j(\mathcal{I})}\ctr_{\underline{\ell},j}^{s_j}=1$.
    A standard computation shows that $z_j \leq t_j$, with equality if and only if $t_j(\underline{\ell}) = t_j$ for all $\underline{\ell}\in\eta_j(\mathcal{I})$.

    Now suppose $d_1<d_2$ and $s_1+t_1>s_2+t_2$ (the opposite case follows analogously).
    Recall from \cref{p:bar-dims} that $\dimB K =\max\{s_1+z_1,s_2+z_2\}$ and (by the assumption above) $\dimA K=s_1+t_1$.
    Therefore, it suffices to show that $z_1<t_1$.

    Assume for contradiction that $z_1=t_1$, so that $t_1(\underline{\ell})=t_1$ for every $\underline{\ell}\in\eta_1(\mathcal{I})$.
    Define $\bm{v}=\bigl(\ctr_{i,1}^{s_1}\ctr_{i,2}^{t_1}\bigr)_{i\in\mathcal{I}}$.
    A direct computation gives that
    \begin{equation*}
        \frac{H(\eta_1(\bm{v}))}{\chi_1(\bm{v})} = s_1\qquad\text{and}\qquad \frac{H(\bm{v})-H(\eta_1(\bm{v}))}{\chi_2(\bm{v})} = t_1.
    \end{equation*}
    If $\chi_1(\bm{v})\leq\chi_2(\bm{v})$ (so that $\bm{v}\in\mathcal{P}_1$), then
    \begin{equation*}
      \dimH K = d_2 > d_1\geq\frac{H(\eta_1(\bm{v}))}{\chi_1(\bm{v})}+\frac{H(\bm{v})-H(\eta_1(\bm{v}))}{\chi_2(\bm{v})}=s_1+t_1 = \dimA K
    \end{equation*}
    which is a contradiction.

    Otherwise, $\chi_1(\bm{v})>\chi_2(\bm{v})$, so that $\bm{v}\in\mathcal{P}_2$.
    We first show that entropy is additive with respect to the projections: $H(\bm{v}) = H(\eta_1(\bm{v})) + H(\eta_2(\bm{v}))$.
    First, since each index $i\in\mathcal{I}$ is uniquely determined by the pair $(\eta_1(i),\eta_2(i))$, entropy is subadditive with respect to the two projections.
    Thus $H(\eta_2(\bm{v}))\geq H(\bm{v})-H(\eta_1(\bm{v}))=t_1\chi_2(\bm{v})$ and therefore
    \begin{align*}
      d_2&\geq\frac{H(\eta_2(\bm{v}))}{\chi_2(\bm{v})}+\frac{H(\bm{v})-H(\eta_2(\bm{v}))}{\chi_1(\bm{v})}\\
      &\geq\frac{H(\bm{v})}{\chi_1(\bm{v})}+t_1\chi_2(\bm{v})\biggl(\frac{1}{\chi_2(\bm{v})}-\frac{1}{\chi_1(\bm{v})}\biggr)\\
      &= s_1+t_1,
    \end{align*}
    where the last inequality uses $H(\bm{v})=s_1\chi_1(\bm{v})+t_1\chi_2(\bm{v})$.
    Since $d_2=\dimH K\leq\dimA K=s_1+t_1$, equality holds throughout.
    Moreover, since $\chi_1(\bm{v})>\chi_2(\bm{v})$, necessarily $H(\eta_2(\bm{v}))=t_1\chi_2(\bm{v})$ and therefore the claim follows.

    Since $\bm{v}$ has full support and again recalling that each index $i\in\mathcal{I}$ is uniquely determined by the pair $(\eta_1(i),\eta_2(i))$, additivity of entropy implies that each equivalence class $\eta_j^{-1}(\underline{\ell})$ must contain exactly one element from each equivalence class $\eta_{j'}^{-1}(\underline{\ell}')$ (see, e.g., \cite[Theorem~2.6.6]{zbl:1140.94001}).
    In words, the carpet is a ``product'' carpet.
    In this case, a standard computation implies that $s_1 + t_1 = s_2 + t_2$, contradicting the original assumption.
\end{proof}
We now prove our main theorem concerning Barański carpets with few large tangents.
\begin{theorem}\label{t:bar-const}
    Let $K$ be a Barański carpet such that $\eta_j(K)$ satisfies the SSC and $\Omega_j\neq\varnothing$ for $j=1,2$.
    Suppose for one of $j=1,2$, $d_j<d_{j'}$ and $\dimB\eta_j(K)+t_j>\dimB\eta_{j'}(K)+t_{j'}$.
    Then
    \begin{equation*}
        \dimH\{x\in K:\dimA(K,x)=\dimA K\}<\dimH K.
    \end{equation*}
\end{theorem}
\begin{proof}
    Suppose $d_1<d_2$ and $\dimB\eta_1(K)+t_1>\dimB\eta_2(K)+t_2$ (the opposite case follows analogously).
    By \cref{p:bar-dims}, $\dimH K=d_2$ and $\dimA K=\dimB\eta_1(K)+t_1$, and \cref{lem:bar-box} gives $\dimB K<\dimA K$.
    Now if $\gamma\in\Omega_2$ and $x=\pi(\gamma)$, then \cref{p:baranski-tans} gives
    \begin{align*}
        \dimA(K,x)&=\max\{\dimB K,\dimB\eta_2(K)+\dimA K_{\eta_2(\gamma),2}\}\\
                  &\leq \max\{\dimB K,\dimB\eta_2(K)+t_2\}<\dimA K.
    \end{align*}
    In particular, if $\dimA(K,x)=\dimA K$, then necessarily $x=\pi(\gamma)$ where $\gamma\in\Omega_0\cup\Omega_1$.
    But $\dimH\pi(\Omega_0\cup\Omega_1)\leq d_1<d_2=\dimH K$, as required.
\end{proof}

The proof of \cref{t:bar-const} applies \cref{p:baranski-tans} only to sequences in $\Omega_{j'}$, where $j$ is the index given in the statement, so that it in fact suffices to assume that $\eta_{j'}(K)$ satisfies the strong separation condition.
As in the Gatzouras--Lalley case, the role of this assumption is to guarantee that every point of $K$ is regular in the sense of the direction-$j'$ analogue of \cref{d:regular}.

\begin{remark}
    In the context of \cref{t:bar-const}, one can in fact prove that the following are equivalent:
    \begin{enumerate}[nl,r]
        \item\label{im:pointwise-small} $\dimH\{x\in K:\dimA(K,x)=\dimA K\}<\dimH K$.
        \item\label{im:tangent-small} $\dimH\{x\in K:\exists F\in\Tan(K,x)\text{ such that }\dimH F=\dimA K\}<\dimH K$.
        \item\label{im:split-directions} For one of $j=1,2$, $d_j<d_{j'}$ and $\dimB\eta_j(K)+t_j>\dimB\eta_{j'}(K)+t_{j'}$.
    \end{enumerate}
    Such a proof follows similarly to the Gatzouras--Lalley case with appropriate modifications to restrict attention only to the family $\Omega_1$ or $\Omega_2$.
    The only additional observation required is that \cite[Lemma~4.3]{zbl:1206.28011} also holds in the Barański case and the uniform subsystem can be chosen so the maps are contracting strictly in direction $j$ and the dimension of the corresponding attractor is arbitrarily close to $d_j$.

    In particular, if one of the above equivalent conditions hold and without loss of generality $d_1>d_2$ and $\dimB\eta_1(K)+t_1<\dimB\eta_2(K)+t_2$, then the Hausdorff dimension of the level set $\varphi(\alpha)=\dimH\{x\in K:\dimA(K,x)=\alpha\}$ is given by the piecewise formula
    \begin{equation*}
        \varphi(\alpha)=\begin{cases}
            \dimH K &: \dimB K\leq \alpha\leq \dimB\eta_1(K)+t_1\\
            d_2 &: \dimB\eta_1(K)+t_1<\alpha\leq \dimA K.
        \end{cases}
    \end{equation*}
    We leave the missing details to the (highly motivated) reader.
\end{remark}
With \cref{t:bar-const} in hand, we can now give an explicit example of a Barański carpet which has few large tangents.
\begin{corollary}\label{c:exceptional-bar}
    There is a Barański carpet $K$ such that
    \begin{equation*}
        \dimH\{x\in K:\dimA(K,x)=\dimA K\}<\dimH K.
    \end{equation*}
\end{corollary}
\begin{proof}
    Fix some $\delta\in [0,1/6)$ and define parameters $\beta=1/4-\delta$, $\alpha_1=1/3-\delta$, and $\alpha_2=1/6-\delta$.
    Now define the families of maps
    \begin{align*}
        \Phi_1 &=\{(x,y)\mapsto (\alpha_1 x, \beta y+i\beta):i=0,\ldots,3\}\\
        \Phi_{2,a}&=\{(x,y)\mapsto (\alpha_2 x+\alpha_1+j \alpha_2, \beta y+i\beta):j=0,1\text{ and }i=0,1\}\\
        \Phi_{2,b}&=\{(x,y)\mapsto (\alpha_2 x+\alpha_1+j \alpha_2, \beta y+i\beta):j=3,4\text{ and }i=2,3\}
    \end{align*}
    and then set
    \begin{equation*}
        \Phi_2=\Phi_{2,a}\cup\Phi_{2,b}\qquad\text{and}\qquad\Phi=\Phi_1\cup\Phi_{2,a}\cup\Phi_{2,b}.
    \end{equation*}
    We abuse notation and use functions and indices interchangeably.
    Note that $\Phi$ is a Barański IFS with five columns; the carpet has the same combinatorial pattern as the carpet generated by the maps depicted in \cref{sf:bar}.
    Both projected IFSs satisfy the SSC whenever $\delta > 1/36$.

    We now simplify the dimensional expression in \cref{p:bar-dims}~\cref{im:bar-adim} for our specific system.
    First, for $\bm{w}\in\mathcal{P}$, set $p=\sum_{i\in\Phi_2}\bm{w}_i$.
    Note that $\chi_1(\bm{w})=-p\log \alpha_2-(1-p)\log\alpha_1$ and $\chi_2(\bm{w})=-\log\beta$ depend only on $p$.
    But since entropy and projected entropy are maximized uniquely by uniform vectors, defining the vector $\bm{z}(p)\in\mathcal{P}$ by
    \begin{equation*}
        \bm{z}(p)_i=\begin{cases}
            \frac{1-p}{4}:i\in\Phi_1\\
            \frac{p}{8}:i\in\Phi_2
        \end{cases}
    \end{equation*}
    we necessarily have
    \begin{align*}
        \frac{H(\eta_1(\bm{w}))}{\chi_1(\bm{w})}+\frac{H(\bm{w})-H(\eta_1(\bm{w}))}{\chi_{2}(\bm{w})} \leq{}&\frac{H(\eta_1(\bm{z}(p)))}{\chi_1(\bm{z}(p))}+\frac{H(\bm{z}(p))-H(\eta_1(\bm{z}(p)))}{\chi_{2}(\bm{z}(p))}\\
        ={}& \frac{-p\log p - (1-p)\log(1-p)+p\log 4}{-p\log\alpha_2 - (1-p)\log\alpha_1}\\
           &+\frac{(2-p)\log 2}{-\log\beta}\\
           &\eqqcolon D_1(p)
    \end{align*}
    and
    \begin{align*}
        \frac{H(\eta_2(\bm{w}))}{\chi_2(\bm{w})}+\frac{H(\bm{w})-H(\eta_2(\bm{w}))}{\chi_{1}(\bm{w})} &\leq\frac{H(\eta_2(\bm{z}(p)))}{\chi_2(\bm{z}(p))}+\frac{H(\bm{z}(p))-H(\eta_2(\bm{z}(p)))}{\chi_{1}(\bm{z}(p))}\\
        &= \frac{-p\log p-(1-p)\log(1-p)+p\log 2}{-p\log \alpha_2-(1-p)\log \alpha_1} + \frac{\log 4}{-\log\beta}\\
        &\eqqcolon D_2(p).
    \end{align*}
    Moreover, writing $p_0=\frac{\log\alpha_1-\log\beta}{\log\alpha_1-\log\alpha_2}$, $\bm{z}(p)\in\mathcal{P}_1$ if and only if $p\in[0,p_0]$ and $\bm{z}(p)\in\mathcal{P}_2$ if and only if $p\in[p_0,1]$.
    We thus observe that
    \begin{equation*}
        \dimH K=\sup_{p\in[0,1]} D(p)\qquad\text{where}\qquad D(p)=\begin{cases}D_1(p) &: 0\leq p\leq p_0\\ D_2(p) &: p_0\leq p\leq 1\end{cases}.
    \end{equation*}
    Now, a direct computation shows that, for $\delta=1/20$, we have $\sup_{p\in[0,1]}D_1(p)=1.45949\ldots$ and $\sup_{p\in[0,1]}D_2(p)=1.47289\ldots$, and moreover the maximum of $D_2$ is attained at the value $p_2=0.53755\ldots$, which lies in $(p_0,1)$ because $p_0=0.39254\ldots$.
    Consequently
    \begin{equation*}
        d_1\leq\sup_{p\in[0,1]}D_1(p)<\sup_{p\in[0,1]} D(p)=D_2(p_2)=d_2.
    \end{equation*}
    (It seems that the same conclusion holds for every $\delta\in(0,1/6)$, but this is not required for the proof.)

    On the other hand, $\dimB\eta_1(K)$ is the unique solution $s$ of $\alpha_1^s+4\alpha_2^s=1$ and the column $\Phi_1$ maximizes $t_1(\underline{\ell})$, so that $t_1=\log 4/\log(1/\beta)$, whereas $\dimB\eta_2(K)=\log 4/\log(1/\beta)$ and $t_2$ is the unique solution $t$ of $\alpha_1^t+2\alpha_2^t=1$.
    For $\delta=1/20$ this gives $\dimB\eta_1(K)+t_1=1.70394\ldots$ and $\dimB\eta_2(K)+t_2=1.47289\ldots$.
    Since $1/20>1/36$ both projected IFSs satisfy the SSC, so that all the conditions of \cref{t:bar-const} are satisfied, as required.
\end{proof}

\begin{acknowledgements}
    The authors thank Balázs Bárány and Lars Olsen for valuable comments on a draft version of this document.
    We also thank the referee for a careful reading of the paper and for helpful suggestions which have improved the exposition.

    AR is supported by Tuomas Orponen's grant from the Research Council of Finland via the project Approximate incidence geometry, grant no.\ 355453.
    Part of this work was also done while AR was a PhD student at the University of St Andrews, supported by EPSRC Grant EP/V520123/1 and the Natural Sciences and Engineering Research Council of Canada.
\end{acknowledgements}

@preprint{arxiv:2107.00983,
  arxiv = {2107.00983},
  author = {Bárány, Balázs and Käenmäki, Antti and Yu, Han},
  eprint = {2107.00983},
  eprinttype = {arxiv},
  journal = {Adv. Math.},
  month = {07},
  title = {Finer geometry of planar self-affine sets},
  year = {2021},
}

@preprint{arxiv:2309.11971,
  arxiv = {2309.11971},
  author = {Käenmäki, Antti and Rutar, Alex},
  eprint = {2309.11971},
  eprinttype = {arxiv},
  month = {09},
  title = {Tangents of invariant sets},
  year = {2023},
}

@article{doi:10.1093/imrn/rnw336,
  author = {Käenmäki, Antti and Ojala, Tuomo and Rossi, Eino},
  doi = {10.1093/imrn/rnw336},
  eprint = {10.1093/imrn/rnw336},
  eprinttype = {doi},
  journal = {Int. Math. Res. Not.},
  month = {6},
  pages = {3769--3799},
  title = {Rigidity of quasisymmetric mappings on self-affine carpets},
  volume = {2018},
  year = {2018},
}

@article{zbl:0757.28011,
  author = {Lalley, Steven P. and Gatzouras, Dimitrios},
  doi = {10.1512/iumj.1992.41.41031},
  eprint = {0757.28011},
  eprinttype = {zbl},
  journal = {Indiana Univ. Math. J.},
  language = {English},
  pages = {533--568},
  title = {Hausdorff and box dimensions of certain self-affine fractals},
  volume = {41},
  year = {1992},
  zbl = {0757.28011},
  zbmath = {00108960},
}

@thesis{zbl:0396.46035,
  address = {Orsay},
  author = {Assouad, Patrice},
  eprint = {0396.46035},
  eprinttype = {zbl},
  institution = {Univ. Paris XI},
  language = {French},
  series = {Publ. Math. Orsay},
  title = {Espaces métriques, plongements, facteurs},
  type = {Thèse de doctorat d’État},
  year = {1977},
  zbl = {0396.46035},
  zbmath = {03615337},
}

@book{zbl:1140.94001,
  address = {Hoboken, NJ},
  author = {Cover, Thomas M. and Thomas, Joy A.},
  doi = {10.1002/047174882X},
  eprint = {1140.94001},
  eprinttype = {zbl},
  language = {English},
  publisher = {John Wiley \& Sons},
  title = {Elements of information theory},
  year = {2006},
  zbl = {1140.94001},
  zbmath = {05060438},
}

@article{zbl:1116.28008,
  author = {Barański, Krzysztof},
  doi = {10.1016/j.aim.2006.06.005},
  eprint = {1116.28008},
  eprinttype = {zbl},
  journal = {Adv. Math.},
  language = {English},
  pages = {215--245},
  title = {Hausdorff dimension of the limit sets of some planar geometric constructions},
  volume = {210},
  year = {2007},
  zbl = {1116.28008},
  zbmath = {05132554},
}

@article{zbl:1154.37322,
  author = {Furstenberg, Hillel},
  doi = {10.1017/S0143385708000084},
  eprint = {1154.37322},
  eprinttype = {zbl},
  journal = {Ergodic Theory Dyn. Syst.},
  language = {English},
  pages = {405--422},
  title = {Ergodic fractal measures and dimension conservation},
  volume = {28},
  year = {2008},
  zbl = {1154.37322},
  zbmath = {05272907},
}

@article{zbl:1206.28011,
  arxiv = {0903.2216v3},
  author = {Ferguson, Andrew and Jordan, Thomas and Shmerkin, Pablo},
  doi = {10.4064/fm209-3-1},
  eprint = {1206.28011},
  eprinttype = {zbl},
  journal = {Fundam. Math.},
  language = {English},
  month = {10},
  origdate = {2009-03-12},
  pages = {193--213},
  title = {The {Hausdorff} dimension of the projections of self-affine carpets},
  version = {3},
  volume = {209},
  year = {2010},
  zbl = {1206.28011},
  zbmath = {05834933},
}

@article{zbl:1278.37032,
  author = {Mackay, John M.},
  doi = {10.1090/S1088-4173-2011-00232-3},
  eprint = {1278.37032},
  eprinttype = {zbl},
  journal = {Conform. Geom. Dyn.},
  language = {English},
  pages = {177--187},
  title = {Assouad dimension of self-affine carpets},
  volume = {15},
  year = {2011},
  zbl = {1278.37032},
  zbmath = {06040469},
}

@article{zbl:1305.28021,
  author = {Fraser, Jonathan M.},
  doi = {10.1090/S0002-9947-2014-06202-8},
  eprint = {1305.28021},
  eprinttype = {zbl},
  journal = {Trans. Am. Math. Soc.},
  language = {English},
  pages = {6687--6733},
  title = {Assouad type dimensions and homogeneity of fractals},
  volume = {366},
  year = {2014},
  zbl = {1305.28021},
  zbmath = {06371220},
}

@article{zbl:1321.54059,
  author = {Le Donne, Enrico and Rajala, Tapio},
  doi = {10.1512/iumj.2015.64.5469},
  eprint = {1321.54059},
  eprinttype = {zbl},
  journal = {Indiana Univ. Math. J.},
  language = {English},
  pages = {21--54},
  title = {Assouad dimension, {Nagata} dimension, and uniformly close metric tangents},
  volume = {64},
  year = {2015},
  zbl = {1321.54059},
  zbmath = {06434506},
}

@book{zbl:1390.28012,
  address = {Cambridge},
  author = {Bishop, Christopher J. and Peres, Yuval},
  doi = {10.1017/9781316460238},
  eprint = {1390.28012},
  eprinttype = {zbl},
  language = {English},
  number = {162},
  publisher = {Cambridge University Press},
  series = {Camb. Stud. Adv. Math.},
  title = {Fractals in probability and analysis},
  year = {2017},
  zbl = {1390.28012},
  zbmath = {06653783},
}

@article{zbl:1364.28011,
  author = {Li, Wenwen and Li, Wenxia and Miao, Junjie and Xi, Lifeng},
  doi = {10.1007/s11464-016-0539-6},
  eprint = {1364.28011},
  eprinttype = {zbl},
  journal = {Front. Math. China},
  language = {English},
  pages = {705--722},
  title = {Assouad dimensions of {Moran} sets and {Cantor}-like sets},
  volume = {11},
  year = {2016},
  zbl = {1364.28011},
  zbmath = {06714858},
}

@book{zbl:1467.28001,
  address = {Cambridge},
  author = {Fraser, Jonathan M.},
  doi = {10.1017/9781108778459},
  eprint = {1467.28001},
  eprinttype = {zbl},
  language = {English},
  number = {222},
  publisher = {Cambridge University Press},
  series = {Camb. Tracts Math.},
  title = {Assouad dimension and fractal geometry},
  year = {2020},
  zbl = {1467.28001},
  zbmath = {07274738},
}

@book{zbl:0819.28004,
  author = {Mattila, Pertti},
  eprint = {0819.28004},
  eprinttype = {zbl},
  language = {English},
  publisher = {Cambridge: Univ. Press},
  series = {Camb. Stud. Adv. Math.},
  subtitle = {Fractals and rectifiability},
  title = {Geometry of sets and measures in {Euclidean} spaces},
  volume = {44},
  year = {1995},
  zbl = {0819.28004},
  zbmath = {00739280},
}

@article{zbl:1527.28004,
  arxiv = {2203.15301},
  author = {Anttila, Roope},
  doi = {10.1017/prm.2022.83},
  eprint = {1527.28004},
  eprinttype = {zbl},
  journal = {Proc. R. Soc. Edinb., Sect. A, Math.},
  language = {English},
  pages = {2053--2078},
  publisher = {Cambridge University Press},
  title = {Pointwise Assouad dimension for measures},
  volume = {153},
  year = {2023},
  zbl = {1527.28004},
  zbmath = {07768297},
}

@article{zbl:1561.28063,
  arxiv = {2209.13952v3},
  author = {Fraser, Jonathan M. and Rutar, Alex},
  doi = {10.54330/afm.142535},
  eprint = {1561.28063},
  eprinttype = {zbl},
  journal = {Ann. Fenn. Math.},
  language = {English},
  month = {12},
  origdate = {2022-09-28},
  pages = {3--21},
  publisher = {Finnish Mathematical Society},
  title = {Assouad-type dimensions of overlapping self-affine sets},
  version = {3},
  volume = {49},
  year = {2024},
  zbl = {1561.28063},
  zbmath = {07808129},
}

@article{zbmath:08103993,
  arxiv = {2410.17944},
  author = {Käenmäki, Antti and Rutar, Alex},
  doi = {10.1017/S0305004125000416},
  eprint = {08103993},
  eprinttype = {zbmath},
  journal = {Math. Proc. Camb. Philos. Soc.},
  language = {English},
  month = {10},
  pages = {623--648},
  publisher = {Cambridge University Press, Cambridge},
  title = {Regularity of non-autonomous self-similar sets},
  volume = {179},
  year = {2025},
  zbmath = {08103993},
}
\end{document}